\documentclass{amsart}  

\usepackage{amsmath,amsfonts,amssymb,amsthm,latexsym}
\usepackage[T1]{fontenc}
\usepackage[latin1]{inputenc}

\usepackage[USenglish]{babel}
\usepackage{mathtools}
\usepackage{graphicx} 
\usepackage{epstopdf}

\makeindex

\usepackage{tikz}
\usetikzlibrary{matrix,arrows,cd}

\usepackage{hyperref}
\usepackage[all]{xy}
\usepackage{url}

\newtheorem{thm}{Theorem}[section]
\newtheorem{prop}[thm]{Proposition}
\newtheorem{lem}[thm]{Lemma}

\theoremstyle{definition}
\newtheorem{defi}[thm]{Definition} 
\newtheorem{example}[thm]{Example} 

\newtheorem{remark}[thm]{Remark}

\numberwithin{equation}{section}

\DeclareMathOperator{\Id}{Id}
\DeclareMathOperator{\car}{char}
 
 \DeclareMathOperator{\Conj}{Conj}
  \DeclareMathOperator{\Ob}{Ob}
   
    \DeclareMathOperator{\Lie}{Lie}

\begin{document}

	\title[Braiding II]{Braiding for categorical algebras and crossed modules of algebras II:\\ Leibniz algebras}

	\author{A. Fern\'andez-Fari\~na}
\address{[A. Fern\'andez-Fari\~na] Department of Matem\'aticas, University of Santiago de Compostela, 15782, Spain.}
\email{alejandrofernandez.farina@usc.es}

	\author{M. Ladra}
\address{[M. Ladra] Department of Matem\'aticas, Institute of Matem\'aticas, University of Santiago de Compostela, 15782, Spain.}
\email{manuel.ladra@usc.es}

\thanks{This work was partially supported by Agencia Estatal de Investigaci\'on (Spain), grant MTM2016-79661-P (European FEDER
support included, UE). The first author is also supported by a scholarship of Xunta de Galicia (Spain), grant ED481A-2017/064, Xunta de Galicia (Spain).}

\begin{abstract}
In this paper we study the category of braided categorical Leibniz algebras and braided crossed modules of Leibniz algebras
and we relate these structures with the categories of braided categorical Lie algebras and braided crossed modules of Lie algebras using the Loday-Pirashvili category.

\end{abstract}
\maketitle
\section*{Introduction}\addcontentsline{toc}{section}{Introduction}

This paper is the immediate continuation of the article~\cite{FFBraidI}. In the first part we study braiding for crossed modules and internal categories of associative and Lie algebras. In this paper we will talk about braidings for the correspondent structures in the Leibniz algebras case.

Leibniz algebras appear in mathematics as a ``non-antisymmetric'' case of Lie algebras. Bearing this in mind, in this paper we will show how to extend the idea of braiding for crossed modules and internal categories of Lie algebras to the Leibniz setting. After introducing these notions, we will prove the equivalence between braided crossed modules of Leibniz algebras and braided categorical Leibniz algebras and we will show the parallelism between its examples and the ones given for groups, associative algebras and Lie algebras, shown in the first part~\cite{FFBraidI}.

For extending the notion of braiding we will use the Loday-Pirashvili category~\cite{LPLM}, allowing us to see Leibniz algebras as a special case of Lie algebras in the category of linear maps with a certain structure of the braided monoidal category.

This second part is organized as follows. In the preliminaries we will recall some definitions from the first part~\cite{FFBraidI} and some basic definitions and properties about crossed modules of Leibniz algebras, including their relationship with crossed modules of Lie algebras. In Section~\ref{S:braidLeib} we show the internalization of the notion of crossed module with a left Lie action of Lie objects in an arbitrary category. We will also define braidings for crossed modules of Lie objects and categorical Lie objects. We will apply this definition to the Loday-Pirashvili category $\mathcal{LM}_K$ and we will obtain the definition of braiding for crossed modules of Leibniz algebras and categorical Leibniz algebras. In Section~\ref{S:Leibequiv}  we will prove the equivalence between braided categories in the Leibniz algebras case. Finally, in Section~\ref{S:nonabeliantensorLeib} we will show that the non-abelian tensor product of Leibniz algebras gives an example of braided crossed module of Leibniz algebras.

\section{Preliminaries}

\subsection{Crossed Modules of Lie algebras}\hfill

We will recall the  definitions of braiding for categorical Lie $K$-algebras and crossed modules of Lie $K$-algebras and some results on them (see~\cite{FFBraidI}).

\begin{defi}
	Let $\mathcal{C}=(C_1,C_0,s,t,e,k)$ be a categorical Lie $K$-algebra.
	
	A\emph{ braiding on $\mathcal{C}$ } is a $K$-bilinear map $\tau\colon C_0\times C_0\xrightarrow{} C_1$, $(a,b)\mapsto \tau_{a,b}$, that verifies the following properties:
	\begin{equation}\label{LieT1}
	\tau_{a,b}\colon [a,b]\xrightarrow{ }[b,a],  \tag{LieT1}
	\end{equation}
	\begin{equation}\label{LieT2}
	\begin{tikzcd}		{[s(x),s(y)]}\arrow[d,"{\tau_{s(x),s(y)}}"]\arrow[r,"{[x,y]}"]& {[t(x),t(y)]}\arrow[d,"{\tau_{t(x),t(y)}}"]\\
	{[s(y),s(x)]}\arrow[r,"{[y,x]}"]& {[t(y),t(x)]},\tag{LieT2}
	\end{tikzcd}
	\end{equation}
	\begin{align}
	\tau_{[a,b],c}&=\tau_{a,[b,c]}-\tau_{b,[a,c]},\tag{LieT3}\label{LieT3}\\
	\tau_{a,[b,c]}&=\tau_{[a,b],c}-\tau_{[a,c],b},\tag{LieT4}\label{LieT4}
	\end{align}
	for $a,b,c\in C_0$, $x,y\in C_1$.
\end{defi}

\begin{defi}
	Let $\mathcal{X}=(M,N,\cdot,\partial)$ be a crossed module of Lie $K$-algebras.
	
	A \emph{braiding} (or \emph{Peiffer lifting}) on the crossed module $\mathcal{X}$ is a $K$-bilinear map
	$\{-,-\}\colon N\times N\xrightarrow{ } M$ verifying:
	\begin{align}
	\partial\{n,n'\}&=[n,n'],\tag{BLie1}\label{BLie1}\\
	\{\partial m, \partial m' \}&=[m,m'],\tag{BLie2}\label{BLie2}\\
	\{\partial m, n \}&=-n\cdot m, \tag{BLie3}\label{BLie3}\\
	\{n,\partial m \}&=n\cdot m, \tag{BLie4}\label{BLie4}\\
	\{n,[n',n'']\}&=\{[n,n'],n''\}-\{[n,n''],n'\},\tag{BLie5}\label{BLie5}\\
	\{[n,n'],n''\}&=\{n,[n',n'']\}-\{n',[n,n'']\},\tag{BLie6}\label{BLie6}
	\end{align}
	for $m,m'\in M$, $n,n',n''\in N$.
	
	If $\{-,-\}$ is a braiding on $\mathcal{X}$ we will say that \emph{$(M,N,\cdot,\partial,\{-,-\})$ is a braided crossed module of Lie $K$-algebras}.
\end{defi}

\begin{defi}
	A $K$-algebra $(M,[-,-])$ is called a \emph{Leibniz }$K$-\emph{algebra}  if the Leibniz identity is verified, i.e.
	\[[x,[y,z]]=[[x,y],z]-[[x,z],y], \qquad  x,y,z\in M.\]
\end{defi}

Remember, as we saw in  \cite{FFBraidI}, that we are working in internal categories where all the internal morphisms are internal isomorphisms, and we have the following properties.

\begin{lem}\label{Lemacomposition}
	Let $(C_1,C_0,s,t,e,k)$ be an internal (associative, Lie, Leibniz) $K$-algebra or a categorical group whose operation is denoted by ``$+$''. Then the following rule for the composition is true
	\[k((x,y))=x-e(t(x))+y=x-e(s(y))+y, \qquad  (x,y)\in C_1\times_{C_0}C_1.\]	
\end{lem}

\begin{prop}[\cite{FFBraidI}]\label{trenzaanticoncor}
	Let $K$ be a field of $\car(K)\neq2$ and $(C_1,C_0,s,t,e,k)$ a categorical Lie $K$-algebra.
	
	If $\tau\colon C_0\times C_0\xrightarrow{}C_1$ is a $K$-bilinear map verifying \eqref{LieT1} and \eqref{LieT2}, then
	\[
	\tau_{a,[b,c]}=[e(a),\tau_{b,c}] \qquad \text{and} \qquad
	\tau_{[b,c],a}=[\tau_{b,c},e(a)].
	\]
	In particular, by the anticommutativity, we have that $\tau_{a,[b,c]}=-\tau_{[b,c],a}$.
\end{prop}

\subsection{Crossed Modules of Leibniz algebras}\hfill

The definition of crossed modules of Leibniz $K$-algebras, ``non-antisymmetric'' case of Lie $K$-algebras,  was introduced by Loday and Pirashvili in \cite{LodayPira}.

\begin{defi}
	Let $N$ and $M$ be two Leibniz $K$-algebras. A \emph{Leibniz action} of $N$ on $M$ is a pair $\cdot=(\cdot_1,\cdot_2)$ where $\cdot_1\colon N\times M\xrightarrow{}M$ and $\cdot_2\colon M\times N\xrightarrow{}M$ are $K$-bilinear maps and the following properties are verified
	\begin{align}
		n\cdot_1[m,m']&=[n\cdot_1 m,m']-[n\cdot_1 m',m],\tag{ALeib1}\label{ALei1}\\
		[m,n\cdot_1 m']&=[m \cdot_2 n,m']-[m,m']\cdot_2 n,\tag{ALeib2}\label{ALei2}\\
		[m,m'\cdot_2 n]&=[m,m']\cdot_2 n-[m\cdot_2 n,m'],\tag{ALeib3}\label{ALei3}\\
		m\cdot_2[n,n']&=(m\cdot_2 n)\cdot_2 n'-(m\cdot_2 n')\cdot_2 n,\tag{ALeib4}\label{ALei4}\\
		n\cdot_1(m\cdot_2 n')&=(n\cdot_1 m)\cdot_2 n'-[n,n']\cdot_1 m,\tag{ALeib5}\label{ALei5}\\
		n\cdot_1(n'\cdot_1 m)&=[n,n']\cdot_1 m-(n\cdot_1 m)\cdot_2n',\tag{ALeib6}\label{ALei6}
	\end{align}
	with $m,m'\in M$, $n,n'\in N$.
\end{defi}

\begin{remark}
	If we change the notation of $\cdot_1$ and $\cdot_2$ by $[-,-]$ in both  cases, the axioms of the Leibniz actions are all possible rewritings of the Leibniz identity when we choose two elements in $M$ and one in $N$ (the first three) or one in $M$ and two $N$ (the last three).
	
	In particular we have that the pair $([-,-],[-,-])$ where $[-,-]$ is the Leibniz bracket of the Leibniz $K$-algebra $M$ is a Leibniz action of $M$ on itself.
\end{remark}

\begin{defi}
	A \emph{crossed module of Leibniz $K$-algebras} is a $4$-tuple $(M,N,\cdot,\partial)$ where $M$ and $N$ are Leibniz $K$-algebras, $\cdot=(\cdot_1,\cdot_2)$ is a Leibniz action of $N$ on $M$, $\partial\colon M\xrightarrow{}N$ is a Leibniz $K$-homomorphism and the following properties are verified for $m,m'\in M$, $n\in N$:
	
	$\partial$ is an $N$-equivariant Leibniz $K$-homomorphism (we suppose that the action is given by the bracket in $N$), i.e.
\begin{align*}
		\partial(n\cdot_1 m) &=[n,\partial(m)] \ \text{and} \ \partial(m\cdot_2 n)=[\partial(m),n],\\
		\partial(m)\cdot_1 m'& =[m,m']=m\cdot_2\partial(m')   \qquad \text{(Peiffer identity)}.
	\end{align*}
\end{defi}

\begin{example}
If $M$ is a Leibniz $K$-algebra then  $(M,M,([-,-],[-,-]),\Id_M)$ is a crossed module of Leibniz $K$-algebras (see \cite{FFBraidI} for the Lie case).
\end{example}

The next propositions give a relation between  crossed modules of Lie and Leibniz $K$-algebras.

\begin{prop}
	Let $M$ and $N$ be two Lie $K$-algebras. Then,
	$\cdot$ is a Lie action of $N$ on $M$ if and only if $(\cdot, \cdot^-)$ is a Leibniz action of $N$ on $M$, where $\cdot^-\colon M\times N\xrightarrow{}M$ is defined by $m\cdot^-n\coloneqq-n\cdot m$.
	
	That is, the Lie action is a particular case of a Leibniz action when the action is ``anticommutative''.
\end{prop}

\begin{prop}\label{LieM->LeibM}
	Let $M$ and $N$  be Lie K-algebras. Then,
	$(M,N,\cdot,\partial)$ is a crossed module of Lie $K$-algebras if and only if $(M,N,(\cdot,\cdot^{-}),\partial)$ is a crossed module of Leibniz $K$-algebras.
\end{prop}

\begin{defi}
	A \emph{homomorphism of crossed modules of Leibniz $K$-algebras} between $(M,N,\cdot,\partial)$ and $(M',N',*,\partial')$
	is a pair of Leibniz $K$-homomorphisms, $f_1\colon M \xrightarrow{}M'$ and $f_2\colon N\xrightarrow{}N'$ such that, for $n\in N$ and $m\in M$,
	\[
	f_1(n\cdot_1 m)=f_2(n)*_1f_1(m), \quad  f_1(m\cdot_2 n)=f_1(m)*_2 f_2(n), \ \ \text{and} \ \ \partial' \circ f_1=f_2\circ \partial.\\
	\]
\end{defi}

	We will denote by $\textbf{\textit{X}}(\textbf{\textit{LeibAlg}}_K)$ the category of crossed modules of Leibniz $K$-algebras and its homomorphisms.

\begin{remark}
	As in the case of groups and Lie $K$-algebras we have an equivalence between the categories $\textbf{\textit{ICat}}(\textbf{\textit{LeibAlg}}_K)$ and $\textbf{\textit{X}}(\textbf{\textit{LeibAlg}}_K)$. A proof of this can be found in \cite{ThRa}.
	
	It is obvious that $\textbf{\textit{X}}(\textbf{\textit{LieAlg}}_K)$ can be seen as a full subcategory of the category $\textbf{\textit{X}}(\textbf{\textit{LeibAlg}}_K)$ using Proposition~\ref{LieM->LeibM} (we actually have a functorial isomorphism between a full subcategory of $\textbf{\textit{X}}(\textbf{\textit{LeibAlg}}_K)$ and $\textbf{\textit{X}}(\textbf{\textit{LieAlg}}_K)$).
	
	Since the pullbacks in $\textbf{\textit{LieAlg}}_K$ and $\textbf{\textit{LeibAlg}}_K$ are the same
	 it is immediate to show that $\textbf{\textit{ICat}}(\textbf{\textit{LieAlg}}_K)$ is a full subcategory of $\textbf{\textit{ICat}}(\textbf{\textit{LeibAlg}}_K)$.
	
	We said this because we know that the equivalence in the case of Leibniz generalizes the equivalence in the case of Lie. The action in the functors (which was given in \cite{ThRa}) is given by the bracket and then, it is anticommutative when we have Lie $K$-algebras. The only thing that we need to check is that the Leibniz semidirect product generalizes the Lie semidirect product, but this is immediate from definition (since $m\cdot_2 n'=-n'\cdot_1 m$ is the Lie case), as we will see below.
\end{remark}

\begin{defi}
	Let $M$ and $N$ be two Leibniz $K$-algebras and $\cdot$ a Leibniz action of $N$ on $M$. The \emph{semidirect product} is the $K$-vector space $M\times N$ with the  bracket
	\[
	[(m,n),(m',n')]\coloneqq ([m,m']+n\cdot_1 m'+m\cdot_2 n',[n,n']).
	\]
	
	We will denote the semidirect product by $M\rtimes N$.
\end{defi}

\section{Braiding for categorical Leibniz algebras and crossed modules of Leibniz algebras}\label{S:braidLeib}

In this section we will use the idea of Loday and Pirashvili (\cite{LPLM}) that allows us to see the Leibniz $K$-algebras as a particular case of a Lie algebra in  the monoidal category of linear maps $\mathcal{LM}_K$, also known as the Loday-Pirashvili category (\cite{FGL,Ro15}). Using this, we will  try to define the concept of braiding in the case of Leibniz algebras taking advantage of the fact that they will be a particular case of braidings for the corresponding ideas over Lie objects in that category.

First, we will introduce some notation.

Let \textbf{\textit{C}}  be a category with coproducts $\oplus$. If we have $A\xrightarrow{f}C\xleftarrow{g}B$ we denote the unique morphism given by the universal property of coproduct as $f\boxplus g\colon A\oplus B\xrightarrow{}C$.
	
We let the notation $\oplus$ in morphisms for the coproduct bifunctor and $+$ for the addition in an additive category.

The definition of the category $\mathcal{LM}_K$  can be seen in \cite{LPLM}.

\begin{defi}
	The category $\mathcal{LM}_K$ is the monoidal category with the following data:
	\begin{itemize}
		\item As objects we take the $K$-linear maps.
		\item If {\scriptsize\begin{tikzcd}M\arrow[d,"f"]\\ N\end{tikzcd}} and {\scriptsize\begin{tikzcd}L\arrow[d,"g"]\\ H\end{tikzcd}} are two linear maps, a morphism between them is a pair $\alpha=(\alpha_1,\alpha_2)$, $\alpha_1\colon M\xrightarrow{}L$ and $\alpha_2\colon N\xrightarrow{}H$ of  $K$-linear maps such that the following diagram is commutative:
		\begin{center}
			\begin{tikzcd}
				M\arrow[r,"\alpha_1"]\arrow[d,"f"] & L \arrow[d,"g"]\\
					N\arrow[r,"\alpha_2"] & H.
			\end{tikzcd}
		\end{center}
		\item The composition of two morphisms and identity morphism is given in the obvious way.
		\item The tensor product
		 {\scriptsize $\begin{tikzcd}M\arrow[d,"f"]\\ N\end{tikzcd} \otimes \begin{tikzcd}L\arrow[d,"g"]\\ H\end{tikzcd} \coloneqq \begin{tikzcd}(M\otimes H)\oplus (N\otimes L)\arrow[d,"(f\otimes \Id_H)\boxplus (\Id_N\otimes g)"]\\ N\otimes H\end{tikzcd}$}
		 where $\otimes$ between $K$-vector spaces is the usual tensor product and $\oplus$ is the direct sum of vector spaces. In morphisms the tensor product is given by $(f_1,f_2)\otimes(g_1,g_2)=((f_1\otimes g_2)\oplus (f_2\otimes g_1),g_1\otimes g_2)$.
	\end{itemize}

This is a braided monoidal category with the braiding given by the isomorphism
$\mathcal{T}_{f,g}=(\mathcal{T}_{f,g}^1,\mathcal{T}_{f,g}^2)\colon$ \!\!\!\!\!\!\! {\scriptsize  $\begin{tikzcd}M\arrow[d,"f"]\\N\end{tikzcd}\otimes
\begin{tikzcd}L\arrow[d,"g"]\\ H\end{tikzcd} \xrightarrow{} \begin{tikzcd}L\arrow[d,"g"]\\ H\end{tikzcd} \otimes \begin{tikzcd}M\arrow[d,"f"]\\ N\end{tikzcd}  $},
with the $K$-linear isomorphism $\mathcal{T}^2_{f,g}=\mathcal{T}_{N,H}\colon N\otimes H\xrightarrow{}H\otimes N$, the usual braiding for the $K$-vector tensor product, and $\mathcal{T}^1_{f,g}\colon (M\otimes H)\oplus (N\otimes L)\xrightarrow{}(L\otimes N)\oplus(H\otimes M)$ given by $\mathcal{T}^1_{f,g}((m\otimes h)+(n\otimes l))=(l\otimes n)+(h\otimes m)$.
\end{defi}

\begin{remark}
	$\mathcal{LM}_K$ is an additive category with {\scriptsize \begin{tikzcd}\{0\}\arrow[d,"0"]\\\{0\}\end{tikzcd}} as zero object, where we have
{\scriptsize $\begin{tikzcd}M\arrow[d,"f"]\\N\end{tikzcd}\times
 \begin{tikzcd}L\arrow[d,"g"]\\H\end{tikzcd}\coloneqq \begin{tikzcd}M\times L\arrow[d,"f\times g"]\\N\times H\end{tikzcd}$} as product and the usual abelian group structure on morphisms.
\end{remark}

As Loday and Pirashvili did in \cite{LPLM} we can categorify the idea of Lie $K$-algebra and define it in the $\mathcal{LM}_K$ category. The categorification is the following one that we generalize to a semigroupal category particular case.

\begin{defi}
	Let $\mathcal{C}=(\textbf{\textit{C}},\otimes,a,\mathcal{T})$ be a braided semigroupal category where \textbf{\textit{C}} is an additive category.
	
	We say that a pair $(A,\mu)$ with $A\in \Ob(\textbf{\textit{C}})$ and $\mu\colon A\otimes A\xrightarrow{}A$ is a \emph{Lie object} in $\mathcal{C}$ if and only if we have that:
	\begin{align*}
	0 &=\mu\circ \mathcal{T}_{A,A}+\mu,\\
	0 &=\mu\circ (\Id_A\otimes \mu)\circ a_{A,A,A}+\mu\circ(\mu\otimes \Id_A)\circ a^{-1}_{A,A,A}\circ (\Id_A\otimes \mathcal{T}_{A,A})\circ a_{A,A,A}\\
& \ \ -\mu\circ(\mu\otimes \Id_A).
	\end{align*}
	
	We will say that a morphism $f\colon A\xrightarrow{}B$ is a \emph{Lie morphism} between $(A,\mu)$ and $(B,\eta)$ if verifies the following diagram:
	\begin{center}
		\begin{tikzcd}
			A\otimes A\arrow[d,"f\otimes f"]\arrow[r,"\mu"] & A\arrow[d,"f"]\\
			B\otimes B \arrow[r,"\eta"] & B.
		\end{tikzcd}
	\end{center}
	
	The composition in \textbf{\textit{C}} of Lie morphisms and the identity morphism of a Lie object
	 in \textbf{\textit{C}} are Lie morphisms. Using this, we can define the category  \textbf{\textit{Lie}}$(\mathcal{C})$.
\end{defi}

\begin{example}
	If we take in \textbf{\textit{Vect}}$_K$ the usual tensor product we have that, if $(V,\mu)$ is a Lie object, then $\mu\colon V\otimes V\xrightarrow{}V$ can be seen as a $K$-bilinear map, $\mu(a,b)=:[a,b]$, which verifies
\begin{gather*}
[a,b]=-[b,a],\\
[a,[b,c]]+[[a,c],b]-[[a,b],c]=0.
\end{gather*}
	
The first is the anticommutativity and the second one is the Leibniz identity.
	
	It is clear that, if $\car(K)\neq 2$, $\textbf{\textit{Lie}}(\textbf{\textit{Vect}}_K)$ is isomorphic to $\textbf{\textit{LieAlg}}_K$.
\end{example}

Since the generalization is only true for $\car(K)\neq 2$, we will assume it for the
rest of the paper.

We want to explain what are an object and a morphism in \textbf{\textit{Lie}}$(\mathcal{LM}_K)$. For that we have the following notation.

\begin{defi}

Let $M$ and $N$ be Lie $K$-algebras and  $\alpha\colon M\xrightarrow{}N$  a Lie $K$-homomorphism.
	
	Let $(V,\cdot)$ be a right $M$-module and $(W,*)$  a right $N$-module. A $K$-linear map $V\xrightarrow{f}W$ is $(\alpha\colon M\xrightarrow{}N,\cdot,*)$\emph{-equivariant} if we have that
	\[
	f(v\cdot m)=f(v)*\alpha(m), \quad \text{for} \ v\in V, \ m\in M.
	\]
	 When $N=M$ and $\alpha=\Id_M$ we said that $f$ is $(M,\cdot,*)$\emph{-equivariant}.
	
	Let $(V,\cdot)$ be a left $M$-module and $(W,*)$  a left $N$-module. A  $K$-linear map $V\xrightarrow{f}W$ is $(\alpha\colon M\xrightarrow{}N,\cdot,*)$\emph{-equivariant} if we have that
	\[
	f(m\cdot v)=\alpha(m)*f(v), \quad \text{for} \ v\in V, \ m\in M.
	\]
 When $N=M$ and $\alpha=\Id_M$ we said that $f$ is $(M,\cdot,*)$\emph{-equivariant}.

	Let $(V,\cdot)$ be a left $M$-module and $(W,*)$  a right $N$-module. A $K$-linear map $V\xrightarrow{f}W$ is $(\alpha\colon M\xrightarrow{}N,\cdot,*)$\emph{-equivariant} if we have that
	\[
	f(m\cdot v)=-f(v)*\alpha(m), \quad \text{for} \ v\in V, \ m\in M.
	\]
	 When $N=M$ and $\alpha=\Id_M$ we said that $f$ is $(M,\cdot,*)$-\emph{equivariant}.
\end{defi}

\begin{remark}
	It is easy to check that if $(M,[-,-])$ is a Lie $K$-algebra, then $(M,[-,-])$ is a right and left $(M,[-,-])$-module.
	
	Also, if $\cdot$ is a Lie action of $N$ in $M$, we have that $(M,\cdot)$ is a left $N$-module.
\end{remark}

Using this, we can see in \cite{LPLM} that a Lie object in $\mathcal{LM}_K$ is the following data:

\begin{defi}
	A \emph{Lie object} in $\mathcal{LM}_K$ is a triple ({\scriptsize \begin{tikzcd}M\arrow[d,"f"]\\ N\end{tikzcd}},$*^M_N,[-,-]_N$) where
	\begin{itemize}
		\item $(N,[-,-]_N)$ is a Lie $K$-algebra.
		\item $*^M_N\colon M\times N\xrightarrow{}M$ is such that $(M,*^M_N)$ is a $(N,[-,-]_N)$-module.
		\item $f$ is $((N,[-,-]_N),*^M_N,[-,-]_N)$-equivariant.
	\end{itemize}
	As in the case of Lie $K$-algebras, we will denote a Lie object in $\mathcal{LM}_K$ using the $K$-linear map
	on which it is defined when there are no confusion.
\end{defi}

\begin{remark}
	The ``anticommutative'' property of Lie object for $\mathcal{LM}_K$ allows to recover the Lie product $\mu=(\mu_1,\mu_2)$ for {\scriptsize\begin{tikzcd}M\arrow[d,"f"]\\ N\end{tikzcd}} with the maps $\mu_2=[-,-]_N$ and $\mu_1\colon(M\otimes N)\oplus(N\otimes M)\xrightarrow{}M$, with $\mu_1((m\otimes n)+(n'\otimes m'))=m*^M_N n-{m'}*^{M}_N n'$.
\end{remark}

\begin{defi}
	Let {\scriptsize\begin{tikzcd}M\arrow[d,"f"]\\ N\end{tikzcd}}  and {\scriptsize\begin{tikzcd}L\arrow[d,"g"]\\ H\end{tikzcd}} be Lie objects. A \emph{Lie morphism} in $\mathcal{LM}_K$ between them is a $\mathcal{LM}_K$ morphism $(\alpha_1,\alpha_2)$ such that:
	\begin{itemize}
		\item $\alpha_2\colon N\xrightarrow{}H$ is a Lie $K$-homomorphism.
		\item $\alpha_1\colon M\xrightarrow{}L$ is a $(\alpha_2\colon N\xrightarrow{}H,*^{M}_N,*^L_H)$-equivariant map.
	\end{itemize}
\end{defi}

In \cite{LPLM} is shown a way to see the Leibniz $K$-algebras as a particular case of Lie objects in $\mathcal{LM}_K$. We show it in the next example.

\begin{example}
	Let $M$ to be a Leibniz $K$-algebra.
	
	We denote for $I_M$ the ideal generated by elements of the form $[x,x]$ with $x\in M$. It is obvious that the quotient Leibniz $K$-algebra is in fact  a Lie $K$-algebra. We will denote its Lie bracket as $\overline{[-,-]}$, and the elements of the quotient as $\overline{m}$ with $m\in M$.
	
	$\Lie(M)\coloneqq\frac{M}{I_M}$ is known as Lieization (note that if $M$ is really a Lie $K$-algebra, then $\Lie(M)$ trivially naturally isomorphic to $M$), and it is obviously functorial.
	
	With this in mind we take the following Lie object in $\mathcal{LM}_K$:
	
	We take {\scriptsize $\begin{tikzcd}M\arrow[d,"\pi_M"]\\ \Lie(M)\end{tikzcd}$} where $\pi(m)=\overline{m}$ is the natural map. It is a Lie object in $\mathcal{LM}_K$ with the following data:
	\begin{itemize}
		\item $m*^M_{\Lie(M)} \overline{m'}=[m,m']$,		\item $[\overline{m},\overline{m'}]_{\Lie(M)}=\overline{[\overline{m},\overline{m'}]}\coloneqq[m,m']$.
	\end{itemize}
	It is obvious that $\pi$ is $(\Lie(M),m*^M_{\Lie(M)},[\overline{m},\overline{m'}]_{\Lie(M)})$-equivariant.
\end{example}

With that example we have a functor $\Phi\colon\textbf{\textit{LeibAlg}}_K\xrightarrow{}\textbf{\textit{Lie}}(\mathcal{LM}_K)$, that is trivially full. This functor is also an injective functor in objects and morphisms, because, as is shown in \cite{LPLM}, we have $\Psi\colon\textbf{\textit{Lie}}(\mathcal{LM}_K)\xrightarrow{}\textbf{\textit{LeibAlg}}_K$ such that $\Psi\circ\Phi=\Id_{\textbf{\textit{LeibAlg}}_K}$. This functor, on objects, can be seen in the following proposition.

\begin{prop}[\cite{LPLM}]
	Let {\scriptsize\begin{tikzcd}M\arrow[d,"f"]\\ N\end{tikzcd}} be a Lie object in $\mathcal{LM}_K$. Then $(M,[-,-])$, where $[m,m'] \coloneqq m*^{M}_N f(m')$, is a Leibniz $K$-algebra.
\end{prop}

In \cite{FGL} we can see that the previous construction can be extended to crossed modules of Lie algebras in $\mathcal{LM}_K$. They did a crossed module with a right action. In this paper we will define which is a crossed module with a left action, or simply a crossed module of Lie objects.

\begin{defi}
	Let $\mathcal{C}=(\textbf{\textit{C}},\otimes,a,\mathcal{T})$ be a braided semigroupal category where \textbf{\textit{C}} is an additive category.

	If $(A,\mu)$ and $(B,\eta)$ are Lie objects, then a \emph{(left) Lie action} of $(B,\eta)$ on $(A,\mu)$ is a morphism $p\colon B\otimes A\xrightarrow{}A$ such that:
	\begin{align*}
	p\circ (\eta\otimes \Id_A)&=p\circ (\Id_B\otimes p)\circ a_{B,B,A}\circ(\Id_{(B\otimes B)\otimes A}-(\tau_{B,B}\otimes \Id_A)),\\
	p\circ (\Id_B\otimes \mu)\circ a_{B,A,A}&=\mu\circ (p\otimes \Id_A)\circ (\Id_{(B\otimes B)\otimes A}-(a^{-1}_{B,A,A}\circ (\Id_{B}\otimes \tau_{A,A})\circ a_{B,A,A})).
	\end{align*}
	
	We said that $((A,\mu),(B,\eta),p,\partial)$ is a \emph{crossed module of Lie objects} if $p$ is a Lie action of $(B,\eta)$ on $(A,\mu)$ and $\partial\colon (A,\mu)\xrightarrow{}(B,\eta)$ is a Lie morphism such that
	\begin{align*}
	\partial\circ p&=\eta\circ (\Id_B\otimes \partial),\\
	\mu&=p\circ(\partial\otimes \Id_A).
	\end{align*}
	
	A \emph{morphism between two crossed modules of Lie objects} $((A,\mu),(B,\eta),p,\partial)$ and $((C,\nu),(D,\theta),q,\delta)$ is a pair of Lie morphisms $(\alpha,\beta)$, $\alpha\colon (A,\mu)\xrightarrow{}(C,\nu)$ and $\beta\colon (B,\eta)\xrightarrow{}(D,\theta)$, which verifies the following diagrams:
	\begin{center}
		\begin{tikzcd}
			B\otimes A\arrow[d,"\beta\otimes \alpha"]\arrow[r,"p"] & A\arrow[d,"\alpha"]\\
			D\otimes C \arrow[r,"q"] & C,
		\end{tikzcd}\qquad
		\begin{tikzcd}
			A\arrow[d,"\alpha"]\arrow[r,"\partial"] & B\arrow[d,"\beta"]\\
			C \arrow[r,"\delta"] & D.
		\end{tikzcd}
	\end{center}
	
	With the usual composition in $\textbf{\textit{C}}\times\textbf{\textit{C}}$ for pairs of morphisms of Lie morphisms we have
	the category \textbf{\textit{XLie}}$(\mathcal{C})$.
\end{defi}

\begin{example}
	If we take in again\textbf{\textit{Vect}}$_K$ with the usual tensor product, since we assume $\car(K)\neq 2$, we have that $\textbf{\textit{XLie}}(\textbf{\textit{Vect}}_K)$ and $\textbf{\textit{X}}(\textbf{\textit{LieAlg}}_K)$ are isomorphic  categories.
\end{example}

We want to know how $\textbf{\textit{XLie}}(\mathcal{LM}_K)$ is, and it is described in the following definition.

\begin{defi}
	Let {\scriptsize\begin{tikzcd}M\arrow[d,"f"]\\ N\end{tikzcd}} and {\scriptsize\begin{tikzcd}L\arrow[d,"g"]\\ H\end{tikzcd}} be Lie objects in $\mathcal{LM}_K$. A \emph{(left) Lie action} of  {\scriptsize\begin{tikzcd}L\arrow[d,"g"]\\ H\end{tikzcd}} on {\scriptsize\begin{tikzcd}M\arrow[d,"f"]\\ N\end{tikzcd}} in $\mathcal{LM}_K$ is a triple $\bar{\cdot}=(\cdot_1,\cdot_2,\xi_\cdot)$ where
	\begin{itemize}
	\item $\cdot_1\colon H\times M\xrightarrow{}M$ is a $K$-bilinear map such that $(M,\cdot_1)$ is a left $H$-module;
	\item $\cdot_2\colon H\times N\xrightarrow{}N$ is a Lie action of $H$ on $N$;
	\item $\xi_\cdot\colon L\times N\xrightarrow{}M$ is a $K$-bilinear map;
	\end{itemize}
	such that the following properties are verified:
	\begin{itemize}
		\item $\cdot_1$ and $\cdot_2$ are compatible actions with $*^{M}_N$. That is for $h\in H$, $n\in N$, $m\in M$ we have
		\[
		h\cdot_1 (m*^{M}_N n)=(h\cdot_1 m)*^M_N n+m*^M_N(h\cdot_2 n);
		 \]
		 \item $f$ is an $(H,\cdot_1,\cdot_2)$-equivariant map;
		 \item $\xi_\cdot$ verifies, for $l\in L$, $n,n'\in N$, $h\in H$ the following equalities
		 \begin{align*}
		 	f(\xi_\cdot(l,n))&=g(l)\cdot_2 n,\\
		 	\xi_\cdot(l*^{L}_H h,n)&=\xi_\cdot(l,h\cdot_2 n)-h\cdot_1\xi_\cdot(l,n),\\
		 	\xi_\cdot(l,[n,n']_N)&=\xi_\cdot(l,n)*^{M}_N n'-\xi_\cdot(l,n')*^M_N n.
		 \end{align*}
	\end{itemize}
\end{defi}

\begin{remark}
	 An action is, in fact, a pair $\bar{\cdot}=(\bar{\cdot}_1,\bar{\cdot}_2)$, with the two maps
\[\bar{\cdot}_1\colon (L\otimes N)\oplus(H\otimes N)\xrightarrow{}M \qquad \text{and} \qquad  \bar{\cdot}_2\colon H\otimes N\xrightarrow{}N\] verifying the general properties, but we can easily obtain the previous definition taking $\cdot_2\coloneqq\bar{\cdot}_2$ and recovering $\bar{\cdot}_1((l\otimes n)+(h\otimes m))=:\xi_\cdot(l,n)+h\cdot_1 m$.
\end{remark}

\begin{defi}
	A \emph{crossed module of Lie objects in $\mathcal{LM}_K$} is a $4$-tuple ({\scriptsize \begin{tikzcd}M\arrow[d,"f"]\\ N\end{tikzcd},\begin{tikzcd}L\arrow[d,"g"]\\ H\end{tikzcd}}$,\bar{\cdot},\partial$) where {\scriptsize\begin{tikzcd}M\arrow[d,"f"]\\ N\end{tikzcd}} and {\scriptsize\begin{tikzcd}L\arrow[d,"g"]\\ H\end{tikzcd}} are Lie objects in $\mathcal{LM}_K$, $\bar{\cdot}$ is a Lie action of {\scriptsize\begin{tikzcd}L\arrow[d,"g"]\\ H\end{tikzcd}} on {\scriptsize\begin{tikzcd}M\arrow[d,"f"]\\ N\end{tikzcd}} and $\partial=(\partial_1,\partial_2)\colon$ \!\!\!{\scriptsize\begin{tikzcd}M\arrow[d,"f"]\\ N\end{tikzcd}}$ \xrightarrow{} $ {\scriptsize\begin{tikzcd}L\arrow[d,"g"]\\ H\end{tikzcd}} is a Lie morphism in $\mathcal{LM}_K$ such that
	\begin{itemize}
		\item $(N,H,\cdot_2,\partial_2)$ is a crossed module of Lie $K$-algebras;
		\item $\partial_1$ is an $(H,\cdot_1,*^{L}_H)$-equivariant map;
		\item The following equalities are true for $h\in H$, $l\in L$, $m\in M$, $n\in N$
		\begin{align*}
		\partial_1(\xi_\cdot(l,n))&=l*^{L}_N \partial_2(h),\\
		\xi_\cdot(\partial_1(m),n)&=m*^{M}_N n=-\partial_2(n)\cdot_1 m.
		\end{align*}
	\end{itemize}
\end{defi}

\begin{defi}
	Let  ({\scriptsize \begin{tikzcd}M\arrow[d,"f"]\\ N\end{tikzcd},\begin{tikzcd}L\arrow[d,"g"]\\ H\end{tikzcd}}$,\bar{\cdot},\partial$) and  ({\scriptsize \begin{tikzcd}X\arrow[d,"k"]\\ Y\end{tikzcd}, \begin{tikzcd}V\arrow[d,"h"]\\ W\end{tikzcd}},$\bar{\star},\delta$) be crossed modules of Lie objects in $\mathcal{LM}_K$. \emph{A morphism of crossed modules of Lie objects in $\mathcal{LM}_K$} between them is a pair $(\alpha,\beta)$ of Lie morphisms $\alpha=(\alpha_1,\alpha_2)\colon$ \!\!\!{\scriptsize\begin{tikzcd}M\arrow[d,"f"]\\ N\end{tikzcd}}$\xrightarrow{}$ {\scriptsize\begin{tikzcd}X\arrow[d,"k"]\\ Y\end{tikzcd}} and $\beta=$ \ $(\beta_1,\beta_2)\colon$~\!\!\!{\scriptsize\begin{tikzcd}L\arrow[d,"g"]\\ H\end{tikzcd}}$\xrightarrow{}${\scriptsize\begin{tikzcd}V\arrow[d,"h"]\\ W\end{tikzcd}} such that
	\begin{itemize}
		\item $(\alpha_2,\beta_2)\colon (N,H,\cdot_2,\partial_2)\xrightarrow{}(Y,W,\star_2,\delta_2)$ is an homomorphism of crossed modules of Lie $K$-algebras;	
			\item $\alpha_1(\xi_\cdot(l,n))=\xi_{\star}(\beta_1(l),\alpha_2(n))$ for $l\in L$, $n\in N$;
			\item $\alpha_1$ is an $(H\xrightarrow{\beta_2}W,\cdot_1,\star_1)$-equivariant map;
			\item	 $\beta_1\circ\partial_1=\delta_1\circ \alpha_1$.
	\end{itemize}
\end{defi}

As in the case of Leibniz $K$-algebras we want to have a pair of functors between the categories $\textbf{\textit{XLie}}(\mathcal{LM}_K)$  and $\textbf{\textit{XLeibAlg}}_K$ who verify similar properties to the simplest case.

For that we do the following propositions of which we omit their proofs because they are immediate using the definitions. The first is symmetrical to the construction we can see in \cite{FGL} for crossed modules with right actions.

\begin{prop}
	Let $(M,N,(\cdot_1,\cdot_2),\partial)$ be a crossed module of Leibniz $K$-algebras.
	
	Then {\rm ({\scriptsize\begin{tikzcd}M\arrow[d,"\pi_M"]\\ \frac{M}{[M,N]_x}\end{tikzcd}},{\scriptsize\begin{tikzcd}N\arrow[d,"\pi_N"]\\ \Lie(N)\end{tikzcd}}$,\bar{\bar{\cdot}},\overline{\partial}$)} is a crossed module of Lie objects in $\mathcal{LM}_K$, where
	\begin{itemize}
	\item $\frac{M}{[M,N]_x}$ is the Lie $K$-algebra quotient of $M$ by the ideal $[M,N]_x$ whose generators are $[m,m]$ for $m\in M$ and $n\cdot_1 m+m\cdot_2 n$ for $n\in N$, $m\in M$, and we denote the natural map by $\pi_M\colon M\xrightarrow{}\frac{M}{[M,N]_x}$, and the elements of that Lie $K$-algebra as $\overline{m}\in\frac{M}{[M,N]_x}$,
	\item $\bar{\cdot}_1\colon \Lie(N)\times M\xrightarrow{}M$, $(\overline{n},m)\mapsto -m\cdot_2 n$,
	\item $\bar{\cdot}_2\colon \Lie(N)\times \frac{M}{[M,N]_x}\xrightarrow{}\frac{M}{[M,N]_x}$, $(\overline{n},\overline{m})\mapsto \overline{n\cdot_1 m}=\overline{-m\cdot_2 n}$,
	\item $\xi_{ \bar{\cdot}}\colon N\times\frac{M}{[M,N]_x}\xrightarrow{}M$, $(n,\overline{m})\xrightarrow{}n\cdot_1 m$,
	\item $\overline{\partial}_1\colon M\xrightarrow{}N$, $m\mapsto \partial(m)$,
	\item $\overline{\partial}_2\colon \frac{M}{[M,N]_x}\xrightarrow{}\Lie(N)$, $\overline{m}\mapsto \overline{\partial m}$.
	\end{itemize} 	
\end{prop}

\begin{remark}
	We will say that the bottom part $(\frac{M}{[M,N]_x},\Lie(N),\bar{\cdot}_2,\overline{\partial}_2)$ is  the Lieization of a crossed module of Leibniz $K$-algebras. In this way we found a similar relation with the Leibniz and Lie object case, who is that in the bottom part we have the Lieization for this functor.
	
	This Lieization verifies again that applied on a crossed module of Lie $K$-algebras, thought as a crossed module of Leibniz $K$-algebras with the action $(\cdot,\cdot^{-})$, is naturally isomorphic to itself. That occurs because, in the quotient, the second generators are null too:
 \[
	n\cdot_1 m+m\cdot_2 n=n\cdot m+m\cdot^{-}n=n\cdot m-n\cdot m=0.
\]
\end{remark}

\begin{prop}
	Let {\rm ({\scriptsize\begin{tikzcd}M\arrow[d,"f"]\\ N\end{tikzcd}},{\scriptsize\begin{tikzcd}L\arrow[d,"g"]\\ H\end{tikzcd}}$,\bar{\cdot},\partial$)} be a crossed module of Lie objects in $\mathcal{LM}_K$, then $(M,L,(\tilde{\cdot}_1,\tilde{\cdot}_2),\partial_1)$ is a crossed module of Leibniz $K$-algebras, where:
	\begin{itemize}
		\item The Leibniz brackets are: $[m,m']=m*^M_N f(m')$ for $m,m'\in M$ and $[l,l']=l*^L_H g(l')$ for $l,l'\in M$,
		\item $\tilde{\cdot_1}\colon L\times N\xrightarrow{}M$ is defined by $l\tilde{\cdot}_1 m=\xi_\cdot(l,f(m))$ for $l\in L$, $m\in M$,
		\item $\tilde{\cdot}_2\colon M\times L\xrightarrow{}M$ is defined by $m\tilde{\cdot_2}l=-g(l)\cdot_1 m$ for $ l\in L$, $m\in M$.
	\end{itemize}
\end{prop}

Thus we have a functor $X\Phi\colon\textbf{\textit{X}}(\textbf{\textit{LeibAlg}}_K)\xrightarrow{}\textbf{\textit{XLie}}(\mathcal{LM}_K)$, that is trivially full, and a functor $X\Psi\colon\textbf{\textit{XLie}}(\mathcal{LM}_K)\xrightarrow{}\textbf{\textit{X}}(\textbf{\textit{LeibAlg}}_K)$. Again we have that $X\Psi\circ X\Phi=\Id_{\textbf{\textit{X}}(\textbf{\textit{LeibAlg}}_K)}$. For that, the functor $X\Phi$ is a full inclusion functor.

\subsection{Braiding for Crossed modules of Lie objects in $\mathcal{LM}_K$ and crossed modules of Leibniz algebras}\hfill

We want to define the notion of braiding for crossed modules of Leibniz algebras. We will do that using $\mathcal{LM}_K$ category and, symmetrical to the previous, we will use the idea that the braiding for crossed module of Leibniz $K$-algebras must be a particular case of braiding for Lie objects in $\mathcal{LM}_K$, verifying symmetrical properties to the previous ones. For this, we define which is a braiding for a crossed modules of Lie objects.

\begin{defi}
	Let $\mathcal{C}=(\textbf{\textit{C}},\otimes,a,\mathcal{T})$ be a braided semigroupal category where \textbf{\textit{C}} is an additive category. Let $\mathcal{X}=((A,\mu),(B,\eta),p,\partial)$ be a crossed module of Lie objects in $\mathcal{C}$.
	
	A \emph{braiding} (or \emph{Peiffer lifting}) on $\mathcal{X}$ is a morphism $\mathfrak{T}\colon B\otimes B\xrightarrow{}A$ verifying:
	\begin{align*}
	\partial\circ \mathfrak{T}&=\eta,\\
	\mathfrak{T}\circ (\partial\otimes\partial)&=\mu,\\
	-\mathfrak{T}\circ (\partial\otimes \Id_B)&=p\circ \mathcal{T}_{A,B},\\
	\mathfrak{T}\circ (\Id_B\otimes\partial)&=p ,
	\end{align*}
	\begin{align*}
	\mathfrak{T}\circ (\Id_B\otimes \eta)\otimes a_{B,B,B}&=\mathfrak{T}\circ (\eta\otimes\Id_B )\circ (\Id_{(B\otimes B)\otimes B}-(a^{-1}_{B,B,B}\circ (\Id_B\otimes \mathcal{T}_{B,B})\circ a_{B,B,B})),\\
	\mathfrak{T}\circ (\eta\otimes \Id_B)&=\mathfrak{T}\circ (\Id_B\otimes \eta)\circ a_{B,B,B}\circ (\Id_{(B\otimes B)\otimes B}-(\mathcal{T}_{B,B}\otimes \Id_B)).
	\end{align*}
	 $((A,\mu),(B,\eta),p,\partial,\mathfrak{T})$ will be called a \emph{braided crossed module of Lie objects} in $\mathcal{C}$.
	
	We will say that a morphism of crossed modules of Lie objects in the category $\mathcal{C}$, $((A,\mu),(B,\eta),p,\partial)\xrightarrow{(\alpha,\beta)}((C,\nu),(D,\theta),q,\delta)$, is a \emph{morphism of braided crossed modules of Lie object} between $((A,\mu),(B,\eta),p,\partial,\mathfrak{T})$ and $((C,\nu),(D,\theta),q,\delta,\mathfrak{Y})$ if it verifies, in addition, that the following diagram commutes
	\begin{center}
		\begin{tikzcd}
			B\otimes B\arrow[d,"\beta\otimes\beta"]\arrow[r,"\mathfrak{T}"] & A\arrow[d,"\alpha"]\\
			D\otimes D \arrow[r,"\mathfrak{Y}"] & B.
		\end{tikzcd}
	\end{center}
		The same composition and identity as in \textbf{\textit{XLie}}$(\mathcal{C})$ is allowed for braided crossed modules of Lie objects in $\mathcal{C}$ and morphisms of braided crossed modules of Lie objects in $\mathcal{C}$. We denote this new category as \textbf{\textit{BXLie}}$(\mathcal{C})$.
\end{defi}

\begin{example}
	As in the previous cases, if we take again \textbf{\textit{Vect}}$_K$ with the usual tensor product, since we assume $\car(K)\neq 2$, we have that $\textbf{\textit{BXLie}}(\textbf{\textit{Vect}}_K)$ and $\textbf{\textit{BX}}(\textbf{\textit{LieAlg}}_K)$ are isomorphic.
\end{example}

$\textbf{\textit{BXLie}}(\mathcal{LM}_K)$ is described in the following definitions.

\begin{defi}
	Let $\mathcal{X}=$({\scriptsize\begin{tikzcd}M\arrow[d,"f"]\\ N\end{tikzcd}},{\scriptsize\begin{tikzcd}L\arrow[d,"g"]\\ H\end{tikzcd}}$,\bar{\cdot},\partial$) be a crossed module of Lie objects in $\mathcal{LM}_K$.
	
	A \emph{braiding} (or \emph{Peiffer lifting}) for $\mathcal{X}$ is  given by a triple of maps $T_{\{-,-\}}=(\{-,-\}_{LH},\{-,-\}_{HL},\{-,-\}_2)$ where:
	\begin{itemize}
		\item $\{-,-\}_2\colon H\times H\xrightarrow{}N$ is a $K$-bilinear map such that $(N,H,\cdot_2,\partial_2,\{-,-\}_2)$ is a braided crossed module of Lie $K$-algebras.
		\item $\{-,-\}_{LH}\colon L\times H\xrightarrow{}M$ and $\{-,-\}_{HL}\colon H\times L\xrightarrow{}M$ are $K$-bilinear maps, which with $\{-,-\}_2$ verify the following properties for $l\in L$, $h,h'\in H$, $m\in M$, $n\in N$:
	\begin{align*}
	f(\{l,h\}_{LH})&=\{g(l),h\}_2, & f(\{h,l\}_{HL})&=\{h,g(l)\}_2,\\
	\partial_1\{l,h\}_{LH}&=l*^L_H h,&\partial_1\{h,l\}_{HL}&=-l*^L_H h,\\
	\{\partial_1(m),\partial_2(n)\}_{LH}&=m*^{M}_N n, & \{\partial_2(n),\partial_1(m)\}_{HL}&=-m*^M_N n,\\
	\{\partial_1(m),h\}_{LH}&=-h\cdot_1 m, &\{\partial_2(n),l\}_{HL}&=-\xi_{\cdot}(l,n),\\
		\{l,\partial_2(n)\}&=\xi_{\cdot}(l,n), &	\{h,\partial_1(m)\}&=h\cdot_1 m,\\
	\end{align*}
	\vspace{-1.1cm}
		\begin{align*}
	\{l,[h,h']_H\}_{LH}&=\{l*^L_H h,h'\}_{LH}-\{l*^L_H h',h\}_{LH},\\
	\{[h,h']_H,l\}_{HL}&=-\{h,l*^L_H h'\}_{HL}-\{l*^L_H h,h'\}_{LH},\\
	\{l,[h,h']_H\}_{LH}&=\{l*^L_H h,h'\}_{LH}+\{h,l*^L_H h'\}_{HL},\\
	\{[h,h']_H,l\}_{HL}&=-\{h,l*^L_H h'\}_{HL}+\{h',l*^L_H h\}_{HL}.
	\end{align*}
	\end{itemize}
	We will say that ({\scriptsize\begin{tikzcd}M\arrow[d,"f"]\\ N\end{tikzcd}}, {\scriptsize\begin{tikzcd}L\arrow[d,"g"]\\ H\end{tikzcd}}$,\bar{\cdot},\partial,T_{\{-,-\}}$) is a braided crossed module of Lie objects in $\mathcal{LM}_K$.
\end{defi}

\begin{remark}
	A braiding is a pair $T_{\{-,-\}}=(T^1_{\{-,-\}},T^2_{\{-,-\}})$, but for simplicity we denote $T^1_{\{-,-\}}\colon (L\otimes H)\oplus (H\otimes L)\xrightarrow{}M$ with $T^1_{\{-,-\}}((l\otimes h)+(h'\otimes l'))=\{l,h\}_{LH}+\{h',l'\}_{HL}$ and $T^2_{\{-,-\}}(h,h')=\{h,h'\}_2$.
\end{remark}

\begin{defi}
	Let ({\scriptsize\begin{tikzcd}M\arrow[d,"f"]\\ N\end{tikzcd}},{\scriptsize\begin{tikzcd}L\arrow[d,"g"]\\ H\end{tikzcd}}$,\bar{\cdot},\partial,T_{\{-,-\}}$) and ({\scriptsize\begin{tikzcd}X\arrow[d,"k"]\\ Y\end{tikzcd}},{\scriptsize\begin{tikzcd}V\arrow[d,"h"]\\ W\end{tikzcd}}$,\bar{\star},\delta,T_{\{-,-\}'}$) be braided crossed modules of Lie objects in $\mathcal{LM}_K$. A \emph{morphism of braided crossed modules of Lie objects in $\mathcal{LM}_K$}  is a pair $(\alpha,\beta)$ such that ({\scriptsize\begin{tikzcd}M\arrow[d,"f"]\\ N\end{tikzcd}},{\scriptsize\begin{tikzcd}L\arrow[d,"g"]\\ H\end{tikzcd}}$,\bar{\cdot},\partial$)$\xrightarrow{(\alpha,\beta)}$({\scriptsize\begin{tikzcd}X\arrow[d,"k"]\\ Y\end{tikzcd}},{\scriptsize\begin{tikzcd}V\arrow[d,"h"]\\ W\end{tikzcd}}$,\bar{\star},\delta$) is a morphism of crossed modules of Lie objects in $\mathcal{LM}_K$ verifying:
	\begin{itemize}
		\item $(\alpha_2,\beta_2)\colon (N,H,\cdot_2,\partial_2,\{-,-\}_2)\xrightarrow{}(Y,W,\star_2,\delta_2,\{-,-\}'_2)$ is an homomorphism of braided crossed modules of Lie $K$-algebras,	
		\item $\alpha_1(\{l,h\}_{LH})=\{\beta_1(l),\beta_2(h)\}'_{VW}$ for $l\in L$, $h\in H$,
		\item $\alpha_1(\{h,l\}_{HL})=\{\beta_2(h),\beta_1(l)\}'_{WV}$ for $l\in L$, $h\in H$.
	\end{itemize}
\end{defi}

We want to use the definition of braiding on crossed modules of Lie objects in $\mathcal{LM}_K$ to obtain a definition for crossed modules of Leibniz $K$-algebras.
 For that we will take a braiding on (\!\!{\scriptsize\begin{tikzcd}M\arrow[d,"\pi_M"]\\ \frac{M}{[M,N]_x}\end{tikzcd}},{\scriptsize\begin{tikzcd}N\arrow[d,"\pi_N"]\\ \Lie(N)\end{tikzcd}}$,\bar{\bar{\cdot}},\overline{\partial}$). If we try to take one $K$-bilinear map $\{-,-\}$
  to do that we will find problems with the form of defining the correspondent maps because we have that the first properties restrict us adding one more quotient that we would like it was trivial for Lie $K$-algebras,
   or if we take it to be trivial, the rest of properties prevent it from being made for a general case in Leibniz $K$-algebras
   (if we take $\{n,\overline{n'}\}_{N\Lie(N)}=\{n,n'\}=\{\overline{n},n'\}_{\Lie(N)N}$ for example, the third and fourth property leads us to prove that $M$ must be Lie $K$-algebra).

For this, as in the case of the two actions, we will take for braiding two $K$-bilinear maps $\{-,-\},\langle-,-\rangle\colon N\times N\xrightarrow{}M$, and define $\{n,\overline{n'}\}_{N\Lie(N)}=\{n,n'\}$, $\{\overline{n},n'\}_{\Lie(N)N}=-\langle n',n\rangle$ and $\{\overline{n},\overline{n'}\}_2=\overline{\{n,n'\}}=\overline{-\langle n',n\rangle}$, where we can see that we introduce a new quotient in $M$
 (we will go into detail later). With this, seeking
to verify the properties,
 we obtain the following definition.

\begin{defi}
	Let $\mathcal{X}=(M,N,(\cdot_1,\cdot_2),\partial)$ be a crossed module of Leibniz $K$-algebras.
	
	A \emph{braiding} (or \emph{Peiffer lifting}) on  $\mathcal{X}$ is a pair $(\{-,-\},\langle-,-\rangle)$ of $K$-bilinear maps $\{-,-\},\langle-,-\rangle\colon N\times N\xrightarrow{}M$, $(n,n')\mapsto \{n,n'\}$ and $(n,n')\mapsto\langle n,n'\rangle$, verifying:
	\begin{align}
	\partial\{n,n'\}&=[n,n']=\partial\langle n,n'\rangle,\tag{BLeib1}\label{LeiB1}\\
	\{\partial m,\partial m' \}&=[m,m']=\langle \partial m,\partial m'\rangle,\tag{BLeib2}\label{LeiB2}\\
	\{\partial m,n \}&=m\cdot_2 n=\langle \partial m,n\rangle,\tag{BLeib3}\label{LeiB3}\\
	\{n,\partial m\}&=n\cdot_1 m=\langle n,\partial m\rangle,\tag{BLeib4}\label{LeiB4}\\
	\{n,[n',n'']\} &=\lbrace[n,n'],n'' \rbrace\!-\!\lbrace [n,n''],n'\rbrace ,\tag{BLeib5}\label{LeiB5}\\
	\langle n,[n',n'']\rangle&=\{[n,n'],n''\}\!-\!\langle[n,n''],n'\rangle,\tag{BLeib6}\label{LeiB6}\\
	\{n,[n',n'']\}&=\{[n,n'],n''\}\!-\!\langle[n,n''],n'\rangle,\tag{BLeib7}\label{LeiB7}\\
	\langle n,[n',n'']\rangle&=\langle[n,n'],n''\rangle\!-\!\langle[n,n''],n'\rangle,\tag{BLeib8}\label{LeiB8}
	\end{align}
	for all $n,n',n''\in N$, $m,m'\in M$.
	
	In this case, we say that $(M,N,(\cdot_1,\cdot_2),\partial,(\{-,-\},\langle-,-\rangle))$ is a \emph{braided crossed module of Leibniz $K$-algebras}.
\end{defi}

\begin{defi}
	An \emph{homomorphism of braided crossed modules of Leibniz $K$-algebras} $\mathcal{X}_1\xrightarrow{(f_1,f_2)}\mathcal{X}_2$, where we have that $\mathcal{X}_1=(M,N,(\cdot_1,\cdot_2),\partial,(\{-,-\},\langle-,-\rangle))$ and $\mathcal{X}_2=(M',N',(*_1,*_2),\partial',(\{-,-\}',\langle-,-\rangle'))$ are braided crossed modules of Lie $K$-algebras, is an homomorphism of between the correspondent crossed modules of Leibniz $K$-algebras verifying:
	\begin{align}
		f_1(\{n,n'\})&=\{f_2(n),f_2(n')\}',\tag{LeibHB1}\label{LHB1}\\
	f_1(\langle n,n'\rangle)&=\langle f_2(n),f_2(n')\rangle',\tag{LeibHB2}\label{LHB2}
	\end{align}
	for $n,n'\in N$.
	
	We denote the category of  braided crossed modules of Leibniz $K$-algebras and its homomorphisms as $\textbf{\textit{BX}}(\textbf{\textit{LeibAlg}}_K)$.
\end{defi}

We want to know how to introduce the braided crossed modules of Lie $K$-algebras as a particular case. The next two properties answer this question:

\begin{prop}
	Let $(M,N,(\cdot_1,\cdot_2),\partial,(\{-,-\},\langle-,-\rangle))$ be a braided crossed module of Leibniz $K$-algebras.
	
	If for all $n,n'\in N$ it is verified that $\{n,n'\}=-\langle n',n\rangle$, then we have the following properties:
	\begin{itemize}
		\item $m\cdot_2 n=-n\cdot_1 m$.
		\item $(M,N,\cdot_1,\partial,\{-,-\})$ is a braided crossed module of Lie $K$-algebras.
	\end{itemize}
	\begin{proof}
		We will check first that $M$ and $N$ are Lie $K$-algebras.
		
		Starting from $N$, if we use \eqref{LeiB1}, we have that, for all $n,n'\in N$, $\partial\langle n,n'\rangle=[n,n']$. Then, if we use that $\langle n,n'\rangle=-\{n',n\}$ we obtain, again for \eqref{LeiB1}:
		\begin{align*}
		[n,n]=\partial\langle n,n'\rangle=-\partial\{n',n\}=-[n,n'].
		\end{align*}
		We conclude that $N$ is a Lie $K$-algebras because the bracket is antisymmetric and we are working in a field of $\car(K)\neq 2$.
		
		Now we take $m,m'\in M$. By \eqref{LeiB2} we have that $\langle \partial m,\partial m'\rangle=[m,m']$. Using again $\langle \partial m,\partial m'\rangle=-\{\partial m',\partial m\}$ and \eqref{LeiB2} we have
		\begin{align*}
		[m,m']=\langle\partial m,\partial m'\rangle=-\{\partial m',\partial m\}=-[m',m].
		\end{align*}
				
		We will check that $m\cdot_2 n=-n\cdot_1 m$. For that we take $m\in M$, $n\in N$. We have the following equalities
		\begin{align*}
		m\cdot_2 n=\langle\partial m,n\rangle=-\{n,\partial m \}=-n\cdot_1 m,
		\end{align*}
		where we used \eqref{LeiB3} in the first equality and \eqref{LeiB4} in the third.
		
		Now, we know that $(M,N,\cdot_1,\partial)$ is a crossed module of Lie $K$-algebras using Proposition~\ref{LieM->LeibM}.
			
		We will prove the equivalences for the axioms of braiding.
		
		The first equality of properties \eqref{LeiB1}--\eqref{LeiB4}
		coincides, respectively, with \eqref{BLie1}--\eqref{BLie4} (in the case of \eqref{LeiB3} remember that $m\cdot_2 n=-n\cdot_1 m$).
		
		The second identity of \eqref{LeiB1} and \eqref{LeiB2} is immediate because of the anticommutativity of the bracket, while the second \eqref{LeiB3} is equivalent to \eqref{BLie4} and the second equality of \eqref{LeiB4} is to \eqref{BLie3} (again using that $n\cdot_1 m=-m\cdot_2 n$).
		
		It is clear that \eqref{LeiB5} and \eqref{BLie5} are identical, and it is very easy to prove that \eqref{LeiB8} is equivalent to \eqref{BLie6}.
		
		To see the last equivalences we must prove an earlier property, which is verified for both braidings under our assumptions:
		
		If $n,n',n''\in N$, then $\{[n,n'],n''\}=-\{n'',[n,n']\}$.
		
		We will start in the Lie case (we suppose we have an action $\cdot$)
		\begin{align*}
		\{[n,n'],n''\}=\{\partial\{n,n'\},n''\}=-n''\cdot\{n,n'\}=-\{n'',\partial\{n,n'\}\}=-\{n'',[n,n']\},
		\end{align*}
		where we use  \eqref{BLie1}, \eqref{BLie3} and \eqref{BLie4}.
		
		In Leibniz case it is not true in general, because we need $m\cdot_2n=-n\cdot_1 m$.
		\begin{align*}
		\{[n,n'],n''\}=\{\partial\{n,n'\},n''\}=\{n,n'\}\cdot_2n''=-n''\cdot_1\{n,n'\}&=-\{n'',\partial\{n,n'\}\}\\
		&=-\{n'',[n,n']\},
		\end{align*}
		where we use \eqref{LeiB1}, \eqref{LeiB3} and \eqref{LeiB4}.
		
		With this property (whose proof can be made only with the final axioms that only involve the braids, but that method has been chosen for its simplicity) we can prove easily the remaining equivalences. That is that \eqref{LeiB6} is equivalent to \eqref{BLie6},
		and \eqref{LeiB7} is equivalent to \eqref{BLie5}.
		In particular $(M,N,\cdot_1,\partial \{-,-\})$ is a braided crossed module of Lie $K$-algebras.
	\end{proof}
\end{prop}

The next two propositions are immediate and the second one gives the construction which yields us to obtain the definition.
\begin{prop}
	Let $M$ and $N$ be Lie $K$-algebras. Then,
	$(M,N,\cdot,\partial,\{-,-\})$ is a crossed module of Lie $K$-algebras if and only if $(M,N,(\cdot,\cdot^{-}),\partial,(\{-,-\},\{-,-\}^{-}))$ is a crossed module of Leibniz $K$-algebras.
	
	 $\cdot^{-}\colon M\times N\xrightarrow{}N$ and $\{-,-\}^{-}\colon N\times N\xrightarrow{} M$ are defined as $m\cdot^{-}n=-n\cdot m$ and $\{n,n'\}^{-}=-\{n',n\}$.
\end{prop}

\begin{prop}
	Let $(M,N,(\cdot_1,\cdot_2),\partial,(\{-,-\},\langle-,-\rangle))$ be a braided crossed module of Leibniz $K$-algebras.
	
	Then {\rm ({\scriptsize\begin{tikzcd}M\arrow[d,"\pi_M"]\\ \frac{M}{\{M,N\}_x}\end{tikzcd}},{\scriptsize\begin{tikzcd}N\arrow[d,"\pi_N"]\\ \Lie(N)\end{tikzcd}}$,\bar{\bar{\cdot}},\overline{\partial},(\{-,-\}_{N\Lie(N)},\{-,-\}_{\Lie(N)N},\{-,-\}_2)$)} is a braided crossed module of Lie objects in $\mathcal{LM}_K$, where:
	\begin{itemize}
		\item $\frac{M}{\{M,N\}_x}$ is the Lie $K$-algebra quotient of $M$ by the ideal $\{M,N\}_x$ whose generators are $[x,x]$ for $x\in M$ and $n\cdot_1 m+m\cdot_2 n$ for $n\in N$, $m\in M$, and $\{n,n'\}+\langle n',n\rangle$ for $n,n'\in N$, and we denote the natural map by $\pi_M\colon M\xrightarrow{}\frac{M}{\{M,N\}_x}$, and the elements as $\overline{m}\in\frac{M}{\{M,N\}_x}$,
		\item $\bar{\cdot}_1\colon \Lie(N)\times M\xrightarrow{}M$, $(\overline{n},m)\mapsto -m\cdot_2 n$,
		\item $\bar{\cdot_2}\colon \Lie(N)\times \frac{M}{\{M,N\}_x}\xrightarrow{}\frac{M}{\{M,N\}_x}$, $(\overline{n},\overline{m})\mapsto \overline{n\cdot_1 m}=\overline{-m\cdot_2 n}$,
		\item $\xi_{ \bar{\cdot}}\colon N\times\frac{M}{\{M,N\}_x}\xrightarrow{}M$, $(n,\overline{m})\xrightarrow{}n\cdot_1 m$,
		\item $\overline{\partial}_1\colon M\xrightarrow{}N$, $m\mapsto \partial(m)$,
		\item $\overline{\partial}_2\colon \frac{M}{\{M,N\}_x}\xrightarrow{}\Lie(N)$, $\overline{m}\mapsto \overline{\partial m}$,
		\item $\{-,-\}_{N\Lie(N)}\colon N\times \Lie(N)\xrightarrow{}M$, $(n,\overline{n'})\mapsto \{n,n'\}$,
		\item $\{-,-\}_{\Lie(N)N}\colon \Lie(N)\times N\xrightarrow{}M$, $(\overline{n},n')\mapsto -\langle n',n\rangle$,
		\item $\{-,-\}_2\colon \Lie(n)\times \Lie(N)\xrightarrow{}M$, $(\overline{n},\overline{n'})\mapsto \overline{\{n,n'\}}=\overline{-\langle n',n\rangle}$.
	\end{itemize} 	
\end{prop}

\begin{remark}
	As in the previous cases, $(\frac{M}{\{M,N\}_x},\Lie(N),\bar{\cdot_2},\overline{\partial}_2,\{-,-\}_2)$ will be called Lieization, and it is functorial.
	
	If we apply this Lieization on a crossed module of Lie $K$-algebras, thought as a crossed module of Leibniz $K$-algebras with the action $(\cdot,\cdot^{-})$ and the braiding $(\{-,-\},\{-,-\}^{-})$, the third generators are null too:
\[
	\{n,n'\}+\langle n',n\rangle=\{n,n'\}+\{n',n\}^{-}=\{n,n'\}-\{n,n'\}=0.
\]
	For that, again, in the Lie case, we obtain a natural isomorphism to itself after doing the Lieization.
\end{remark}

\begin{prop}
	Let {\rm ({\scriptsize\begin{tikzcd}M\arrow[d,"f"]\\ N\end{tikzcd}},{\scriptsize\begin{tikzcd}L\arrow[d,"g"]\\ H\end{tikzcd}}$,\bar{\cdot},\partial,T_{\{-,-\}}$)} be a braided crossed module of Lie objects in $\mathcal{LM}_K$, then $(M,L,(\tilde{\cdot}_1,\tilde{\cdot}_2),\partial_1,(\{-,-\}_{T_{\{-,-\}}},\langle-,-\rangle_{T_{\{-,-\}}}))$ is a braided crossed module of Leibniz $K$-algebras, where:
	\begin{itemize}
		\item The Leibniz brackets are: $[m,m']=m*^M_N f(m')$ for $m,m'\in M$ and $[l,l']=l*^L_H g(l')$ for $l,l'\in M$;
		\item $\tilde{\cdot_1}\colon L\times N\xrightarrow{}M$ is defined by $l \ \tilde{\cdot}_1 m=\xi_\cdot(l,f(m))$ for $l\in L$, $m\in M$;
		\item $\tilde{\cdot}_2\colon M\times L\xrightarrow{}M$ is defined as $m\tilde{\cdot}_2l=-g(l)\cdot_1 m$ for $ l\in L$, $m\in M$;
		\item $\{-,-\}_{T_{\{-,-\}}}\colon L\times L\xrightarrow{}M$ is defined as $\{l,l'\}_{T_{\{-,-\}}}=\{l,g(l')\}_{LH}$ for $l,l'\in L$;
		\item $\langle-,-\rangle_{T_{\{-,-\}}}\colon L\times L\xrightarrow{}M$ is defined as $\langle l,l'\rangle_{T_{\{-,-\}}}=-\{g(l'),l\}_{HL}$ for $l,l'\in L$.
	\end{itemize}
\end{prop}

Thus, we have again a functor $BX\Phi\colon\textbf{\textit{BX}}(\textbf{\textit{LeibAlg}}_K)\xrightarrow{}\textbf{\textit{BXLie}}(\mathcal{LM}_K)$, that is full, and $BX\Psi\colon\textbf{\textit{BXLie}}(\mathcal{LM}_K)\xrightarrow{}\textbf{\textit{BX}}(\textbf{\textit{LeibAlg}}_K)$, verifying $BX\Psi\circ BX\Phi=\Id_{\textbf{\textit{BX}}(\textbf{\textit{LeibAlg}}_K)}$. That is, again, that the functor $BX\Phi$ is a full inclusion functor.

\subsection{Braiding for categorical Lie objects in $\mathcal{LM}_K$ and categorical Leibniz algebras}\hfill

We want to also define a braiding for categorical Leibniz $K$-algebras. We will use, as in the crossed module case, the idea of the category $\mathcal{LM}_K$. Since the properties of internal categories are given by diagrams, we only need to show that $\mathcal{LM}_K$ is a category with pullbacks.

\begin{remark}
	$\mathcal{LM}_K$ is a category with pullbacks.
	
	 If we have the morphisms {\scriptsize\begin{tikzcd}A\arrow[d,"f"]\\B\end{tikzcd}}
$\xrightarrow{\alpha}${\scriptsize\begin{tikzcd}X\arrow[d,"h"]\\Y\end{tikzcd}}$\xleftarrow{\beta}${\scriptsize\begin{tikzcd}C\arrow[d,"g"]\\D\end{tikzcd}} then ({\scriptsize\begin{tikzcd}A\times_X C\arrow[d,"f\times_h g"]\\B\times_Y D\end{tikzcd}}$,(\pi_A,\pi_B),(\pi_C,\pi_D)$) is their pullback.

It is easy to check that \textbf{\textit{Lie}}$(\mathcal{LM}_K)$ has the same pullback with the operations $[(b,d),(b',d')]_{B\times _Y D}\coloneqq([b,b']_B,[d,d']_D)$ and $(a,c)*^{A\times_X C}_{B\times_Y D}(b,d)\coloneqq(a*^A_B b, c*^C_D d)$.  So, we can speak  about categorical Lie objects in $\mathcal{LM}_K$.
\end{remark}

As in the crossed module case, we have the following results.

\begin{prop}
	If $(C_1,C_0,s,t,e,k)$ is a categorical Leibniz $K$-algebra then
{\rm (\!\!\!{\scriptsize\begin{tikzcd}C_1\arrow[d,"\pi_{C_1}"]\\\Lie(C_1)\end{tikzcd}},\!{\scriptsize\begin{tikzcd}C_0\arrow[d,"\pi_{C_0}"]\\\Lie(C_0)\end{tikzcd}}$, (s,\Lie(s)),(t,\Lie(t)),(e,\Lie(e)),(k,\overline{k})$)} is a categorical Lie object in $\mathcal{LM}_K$, where we denote the Lieization functor as $\Lie\colon \textbf{\textit{LeibAlg}}_K\xrightarrow{}\textbf{\textit{LieAlg}}_K$,  the composition morphism $\overline{k}\colon \Lie(C_1)\times_{\Lie(C_0)}\Lie(C_1)\xrightarrow{}\Lie(C_1)$ is defined as $\overline{k}(\overline{x},\overline{y})=\overline{\mathring{k}(x,y)}$ and $\mathring{k}\colon C_1\times C_1\xrightarrow{}C_1$ is the extension of $k$ to the Leibniz product, defined as $\mathring{k}(x,y)=x+y-e(s(y))$.
\end{prop}

\begin{remark}
	$\mathring{k}$ is an extension of $k$, since the same formula is verified for composition, as can be seen in  Lemma~\ref{Lemacomposition}.
	
	One can ask why not to extend $k$ as $\mathring{k}'(x,y)=x+y-e(t(x))$, which is not identical to $\mathring{k}$ in general case. We can do it, but for the case that makes us aware the result is the same, because we have
	\[
	\overline{\mathring{k}(x,y)}=\overline{x}+\overline{y}-\Lie(e)(\Lie(s)(\overline{y}))=\overline{x}+\overline{y}-\Lie(e)(\Lie(t)(\overline{x}))=\overline{\mathring{k}'(x,y)},
	\]
	since $(\overline{x},\overline{y})\in \Lie(C_1)\times_{\Lie(C_0)}\Lie(C_1)$ implies $\Lie(s)(\overline{y})=\Lie(t)(\overline{x})$.
\end{remark}

\begin{remark}
	As was expected, we again have in the bottom part the Lieization, and in the case of Lie $K$-algebras thought as Leibniz $K$-algebras we obtain something trivially isomorphic to the identity.
\end{remark}

\begin{prop}
	If {\rm ({\scriptsize\begin{tikzcd}C_1\arrow[d,"f_1"]\\D_1\end{tikzcd}}
	,{\scriptsize\begin{tikzcd}C_0\arrow[d,"f_0"]\\D_0\end{tikzcd}}$,s,t,e,k$)} is a categorical Lie object in
 $\mathcal{LM}_K$, then $(C_1,C_0,s_1,t_1,e_1,k_1)$ is a categorical Leibniz $K$-algebra, where $[x,y]_{C_1}=x*^{C_1}_{D_1}y$ and $[a,b]_{C_0}=a*^{C_0}_{D_0}b$ for $x,y\in C_1, a,b\in C_0$.
\end{prop}

As in the crossed modules case, we can make a pair of functors.

One functor $I\Phi\colon\textbf{\textit{ICat}}(\textbf{\textit{LeibAlg}}_K)\xrightarrow{}\textbf{\textit{ICat}}(\textbf{\textit{Lie}}(\mathcal{LM}_K))$ that is full, and another $I\Psi\colon\textbf{\textit{ICat}}(\textbf{\textit{Lie}}(\mathcal{LM}_K))\xrightarrow{}\textbf{\textit{ICat}}(\textbf{\textit{LeibAlg}}_K)$, verifying $I\Psi\circ I\Phi=\Id_{\textbf{\textit{ICat}}(\textbf{\textit{LeibAlg}}_K)}$. That is, as is usual, that the functor $I\Phi$ is a full inclusion functor.

This new inclusion functor allows us to define a braiding on categorical Leibniz $K$-algebras using again the idea of braiding of Lie objects in $\mathcal{LM}_K$. For that, we will prove a proposition about the existence pullbacks in categories of Lie objects.

\begin{prop}
	Let $\mathcal{C}=(\textbf{\textit{C}},\otimes,a,\mathcal{T})$ be a braided semigroupal category where \textbf{\textit{C}} is an additive category with pullbacks.
	
	Then $Lie(\mathcal{C})$ has pullbacks and if we have two Lie morphisms
	\begin{center}
		\begin{tikzcd}
			& (A,\mu_A)\arrow[d,"f"]\\
			(B,\mu_B)\arrow[r,"g"]&(C,\mu_C),
		\end{tikzcd}
	\end{center}
	the pullback is given by $((A\times_C B,\mu_{A\times_C B}),\pi_A,\pi_B)$ where
	$A\times_C B$ is the pullback in \textbf{C}
	\begin{center}
		\begin{tikzcd}
	A\times_C B\arrow[r,"\pi_A"]\arrow[d,"\pi_B"] & A\arrow[d,"f"]\\
   B\arrow[r,"g"]&C,
	 \end{tikzcd}
	\end{center}
 and $\mu_{A\times_C B}$ is the unique morphism such that $\pi_X \circ \mu_{A\times_C B}=\mu_X\circ (\pi_X\otimes \pi_X)$ for $X\in \{A,B\}$ given for the universal property of pullbacks in \textbf{C} in the following diagram:
 \begin{center}
 	\begin{tikzcd}
 		 (A\times_CB)\otimes (A\times_C B) \arrow[dr,dashrightarrow,"\mu_{A\times_C B}"]\arrow[dd,"{\pi_B\otimes \pi_B}"]\arrow[rr,"{\pi_A\otimes \pi_A}"]&& A\otimes A\arrow[d,"\mu_A"]\\
 			&A\times_C B\arrow[r,"\pi_A"]\arrow[d,"\pi_B"] & A\arrow[d,"f"]\\
 		B\otimes B\arrow[r,"\mu_B"]&B\arrow[r,"g"]&C
 	\end{tikzcd}
 \end{center}
\begin{proof}
	The first thing that we must check is that $\mu_{A\times_C B}$ is well defined. For that we need to check that we can use the universal property, that is: $f\circ \mu_A\circ (\pi_A\otimes \pi_A)=g\circ \mu_B\circ (\pi_B\otimes \pi_B)$. For this we will use that $f$ and $g$ are Lie morphisms, that is $\mu_C\circ (f\otimes f)=f\circ \mu_A$, $\mu_C\circ (g\otimes g)=g\circ\mu_B$, the fact that $\otimes$ is a functor and $f\circ \pi_A=g\circ \pi_B$:
	\begin{align*}
		&f\circ \mu_A\circ (\pi_A\otimes \pi_A)=\mu_C\circ (f\otimes f)\circ (\pi_A\otimes \pi_A)=\mu_C\circ (f\circ\pi_A\otimes f\circ\pi_A)\\
		&=\mu_C\circ (g\circ\pi_B\otimes g\circ\pi_B)=\mu_C\circ (g\otimes g)\otimes (\pi_B\otimes\pi_B)
		=g\circ \mu_B\circ (\pi_B\otimes \pi_B).
	\end{align*}
	That is, $\mu_{A\times_C B}$ is well defined. If $(A\times_B C,\mu_{A\times_C B})$ is a Lie object, then for the definition of the morphism, we have that $\pi_X \circ \mu_{A\times_C B}=\mu_X\circ (\pi_X\otimes \pi_X)$ for $X\in \{A,B\}$, that is, $\pi_A$ and $\pi_B$ will be Lie morphisms. For that we need, before proving that the universal property is verified for pullbacks in $Lie(\mathcal{C})$, only to prove that $(A\times_C B,\mu_{A\times_C B})$ is a Lie object. For simplicity of notation, we will denote $D \coloneqq A\times_C B$.
	
	We will use universal properties to prove the two axioms of Lie object.
	
	We will prove the first for $(D,\mu_D)$. If we take $X\in \{A,B\}$, we have:
	\[
	\pi_X\circ (-\mu_D\circ \mathcal{T}_{D,D})=-\pi_X\circ\mu_D\circ \mathcal{T}_{D,D}=-\mu_X\circ (\pi_X\otimes \pi_X)\circ \mathcal{T}_{D,D}.
	\]
	Since $\mathcal{T}$ is a natural isomorphism, we have that the following diagram is commutative:
	\begin{center}
		\begin{tikzcd}
			D\otimes D\arrow[d,"\pi_X\otimes\pi_X"]\arrow[r,"\mathcal{T}_{D,D}"]& D\otimes D\arrow[d,"\pi_X\otimes\pi_X"]\\
			X\otimes X\arrow[r,"\mathcal{T}_{X,X}"]& X\otimes X.
		\end{tikzcd}
	\end{center}
	Using this, and that $(X,\mu_X)$ is a Lie object, we get
	\[
	\pi_X\circ (-\mu_D\circ \mathcal{T}_{D,D})=-\mu_X \circ \mathcal{T}_{X,X}\circ (\pi_X\otimes \pi_X)=\mu_X\circ (\pi_X\otimes \pi_X).
	\]

	Since $\mu_D$ is the unique morphism that verifies the previous equality for $ X\in\{A,B\}$ we conclude that $\mu_D=-\mu_D\circ \mathcal{T}_{D,D}$.
	
	To prove the second one we will denote as ${}_H0_J\colon H\xrightarrow{}J $ the zero morphism between $H$ and $J$. With this notation is easy to check that the zero morphism in the second axiom of Lie objects is ${}_{(Y\otimes Y)\otimes Y}0_Y$ for a Lie object $(Y,\mu_Y)$.
	
	Since we have that $f\circ{}_{(D\otimes D)\otimes D}0_A={}_{(D\otimes D)\otimes D}0_C=g\circ{}_{(D\otimes D)\otimes D}0_B$ the universal property tell us that there is a unique morphism $k\colon (D\otimes D)\otimes D\xrightarrow{}D$ such that $\pi_X\circ k={}_{(D\otimes D)\otimes D}0_X$.
	
	 This unique morphism is, in fact, ${}_{(D\otimes D)\otimes D}0_{D}$ by the properties of zero morphisms. Using this universal property we will prove that $\mu_D$ verifies the second axiom. We will take $X\in \{A,B\}$.
	
	 For the first summand we have:
	 \begin{align*}
	 &\pi_X\circ \mu_D\circ (\Id_D\otimes\mu_D)\circ a_{D,D,D}=\mu_X\circ (\pi_X\otimes \pi_X)\circ(\Id_D\otimes\mu_D)\circ a_{D,D,D}\\
	 &=\mu_X\circ (\pi_X\otimes (\pi_X\circ\mu_D))\circ a_{D,D,D}=\mu_X\circ  (\pi_X\otimes (\mu_X\circ (\pi_X\otimes \pi_X)))\circ a_{D,D,D}\\
	 &=\mu_X\circ  (\Id_X\otimes \mu_X )\circ (\pi_X\otimes(\pi_X\otimes \pi_X))\circ a_{D,D,D}.
	 \end{align*}
	 $a$ is a natural isomorphism. This means that the following diagram is commutative:
	 \begin{center}
	 	\begin{tikzcd}
	 		(D\otimes D)\otimes D\arrow[r,"a_{D,D,D}"]\arrow[d,"(\pi_X\otimes \pi_X)\otimes \pi_X"]& D\otimes (D\otimes D)\arrow[d,"\pi_X\otimes (\pi_X\otimes \pi_X)"]\\
	 		(X\otimes X)\otimes X\arrow[r,"a_{X,X,X}"]& X\otimes(X\otimes X).
 	 	\end{tikzcd}
	 \end{center}
	 Using this we have that
	 \[
	 \pi_X\circ \mu_D\circ (\Id_D\otimes\mu_D)\circ a_{D,D,D}=\mu_X\circ  (\Id_X\otimes \mu_X )\circ a_{X,X,X}\circ ((\pi_X\otimes\pi_X)\otimes \pi_X).
	 \]
	
	Doing the same for the second (the naturalness of $a$ gives the same naturalness to $a^{-1}$) and third summands we have that:
	\begin{align*}
	&\pi_X\circ\mu_D\circ(\mu_D\otimes \Id_D)\circ a^{-1}_{D,D,D}\circ (\Id_D\otimes \mathcal{T}_{D,D})\circ a_{D,D,D}\\&=\mu_X\circ(\mu_X\otimes \Id_X)\circ a^{-1}_{X,X,X}\circ (\Id_X\otimes \mathcal{T}_{X,X})\circ a_{X,X,X}\circ((\pi_X\otimes \pi_X)\otimes\pi_X),
	\end{align*}
	\[
	\pi_X\circ(-\mu_D\circ(\mu_D\otimes \Id_D))=-\mu_X\circ(\mu_X\otimes \Id_X)\circ ((\pi_X\otimes \pi_X)\otimes\pi_X).
	\]
	If we denote for $\mathcal{L}_Y$ the morphism which is in the first term of the equality of the second axiom of Lie objects for a Lie object $(Y,\mu_Y)$, we have that adding the three last equalities we have the equality $\pi_X\circ \mathcal{L}_D=\mathcal{L}_X\circ ((\pi_X\otimes \pi_X)\otimes\pi_X)$, using the distributivity of the composition. Since $(X,\mu_X)$ is a Lie object, we have that $\mathcal{L}_X={}_{(X\otimes X)\otimes X}0_X$. For this we have that $\pi_X\circ \mathcal{L}_D={}_{(D\otimes D)\otimes D}0_X$.
	
	Now, by the universal property, we have that $\mathcal{L}_D={}_{(D\otimes D)\otimes D}0_{D}$ and $(D,\mu_D)$ is a Lie object.
	
	To finish we need to prove the universal property.
	We only need to check that the morphism given by the pullback in \textbf{\textit{C}} is a Lie morphism in our case.
	
	Let $(Y,\mu_Y)\xrightarrow{h_X}(X,\mu_X)$ with $X\in \{A,B\}$ be two Lie morphisms, such that $f\circ h_A=g\circ h_B$. Then we have a unique morphism in \textbf{\textit{C}} $h\colon Y\xrightarrow{}D$ such that $\pi_X\circ h=h_X$. We will show that $h$ is a Lie morphism.
	
	For that we take the morphism $\pi_X\circ h\circ \mu_Y\colon Y\otimes Y\xrightarrow{}X$. Is trivial that $f\circ\pi_A\circ h\circ \mu_Y=g \circ\pi_A\circ h\circ \mu_Y$, which yields us to say that there is a unique morphism $u\colon Y\otimes Y\xrightarrow{}D$ such that $\pi_X\circ u=\pi_X\circ h\circ \mu_Y$ for $X\in \{A,B\}$. We have immediately that $u= h\circ \mu_Y$.
	
	On the other hand we have
	\[
	\pi_X\circ \mu_D\circ (h\otimes h)=\mu_X\circ (\pi_X\otimes \pi_X)\circ (h\otimes h)=\mu_X\circ (h_X\otimes h_X)=h_X\circ \mu_Y	=\pi_X\circ h\circ \mu_Y,
	\]
	where we use that $h_X$ is a Lie morphism for $X\in \{A,B\}$. Since $u$ is the unique who verifies the equalities, we have that $\mu_D\circ (h\otimes h)= h\circ \mu_Y$, and $h$ is a Lie morphism.
\end{proof}
\end{prop}

We define now what is a braided categorical Lie object.

\begin{defi}
	Let $\mathcal{C}=(\textbf{\textit{C}},\otimes,a,\mathcal{T})$ be a braided semigroupal category where \textbf{\textit{C}} is an additive category with pullbacks.
	
	 Let $\mathfrak{C}=((C_1,\mu_{C_1}),(C_0,\mu_{C_0}),s,t,e,k)$ be a categorical Lie object in $Lie(\mathcal{C})$.
	
	A \emph{braiding} on $\mathfrak{C}$ is a morphism $\tau\colon C_0\otimes C_0\xrightarrow{}C_1$ verifying:
	\begin{itemize}
		\item $
		s\circ \tau=\mu_{C_0} \ \text{and} \ t\circ \tau=\mu_{C_0}\circ \mathcal{T}_{C_0,C_0},
		$
		\item We define $C_0\otimes C_0\xrightarrow{\mu_{C_1}\times_{C_0} (\tau\circ (t\otimes t)), (\tau\circ (s\otimes s))\times_{C_{0}} (\mu_{C_1}\circ \mathcal{T})} C_1\times_{C_0} C_1$ as the two unique morphisms which verify the universal property,
		respectively, in the following diagrams:
		\[
				\begin{tikzcd}
				C_1\otimes C_1\arrow[rrd,bend left,"\mu_{C_1}"]\arrow[rdd,bend right,"\tau\circ(t\otimes t)"']\arrow[rd,dashrightarrow]\\
				&  C_1\times_{C_0}C_1\arrow[r,"\pi_1"]\arrow[d,"\pi_2"] & C_1\arrow[d,"t"]\\
				&   C_1\arrow[r,"s"]& C_0.
			\end{tikzcd}
			\begin{tikzcd}
			C_1\otimes C_1\arrow[rrd,bend left,"\tau\circ (s\otimes s)"]\arrow[rdd,bend right,"\mu_{C_1}\circ\mathcal{T}_{C_1,C_1}"']\arrow[rd,dashrightarrow]\\
			&  C_1\times_{C_0}C_1\arrow[r,"\pi_1"]\arrow[d,"\pi_2"] & C_1\arrow[d,"t"]\\
			&   C_1\arrow[r,"s"]& C_0.
		\end{tikzcd}
		\]
		The universal property is allowed using the first property of this definition and the fact that $s$ an $t$ are Lie morphisms. For the second one, we need too the naturalness of $\mathcal{T}$.
		
		Then it must be satisfied that
		\[
		k\circ (\mu_{C_1}\times_{C_0} (\tau\circ (t\otimes t)))=k\circ ((\tau\circ (s\otimes s))\times_{C_{0}} (\mu_{C_1}\circ \mathcal{T})).
		\]
	\item It must verify
	\begin{align*}
	&\tau\circ (\Id_{C_0}\otimes \mu_{C_0})\otimes a_{C_0,C_0,C_0}\\
	&=\tau\circ (\mu_{C_0}\otimes\Id_{C_0} )\circ (\Id_{(C_0\otimes C_0)\otimes C_0}-(a^{-1}_{C_0,C_0,C_0}\circ (\Id_{C_0}\otimes \mathcal{T}_{C_0,C_0})\circ a_{C_0,C_0,C_0})),\\
	&\tau\circ (\mu_{C_0}\otimes \Id_{C_0})=\tau\circ (\Id_{C_0}\otimes \mu_{C_0})\circ a_{C_0,C_0,C_0}\circ (\Id_{(C_0\otimes C_0)\otimes C_0}-(\mathcal{T}_{C_0,C_0}\otimes \Id_{C_0})).
	\end{align*}
	\end{itemize}
	We will say that $((C_1,\mu_{C_1}),(C_0,\mu_{C_0}),s,t,e,k,\tau)$ is a \emph{braided categorical Lie object} in $\mathcal{C}$.
	
	A \emph{braided internal functor of braided categorical Lie objects in $\mathcal{C}$}, denoted by $((C_1,\mu_{C_1}),(C_0,\mu_{C_0}),s,t,e,k,\tau)\xrightarrow{(F_1,F_0)}((C'_1,\mu_{C'_1}),(C'_0,\mu_{C'_0}),s',t',e',k',\tau')$, is an internal functor, between the respective internal objects, which verifies, in addition, the following diagram:
	\begin{center}
		\begin{tikzcd}
			C_0\otimes C_0\arrow[d,"F_0\otimes F_0"]\arrow[r,"\tau"] & C_1\arrow[d,"F_1"]\\
			C'_0\otimes C'_0 \arrow[r,"\tau'"] & C_1.
		\end{tikzcd}
	\end{center}
	The same composition and identity as in $\textbf{\textit{ICat}}(\textbf{\textit{Lie}}(\mathcal{C}))$ is allowed for braided categorical Lie objects in $\mathcal{C}$ and its braided internal functors. We denote this new category as \textbf{\textit{BICat}}$(\textbf{\textit{Lie}}(\mathcal{C}))$.
\end{defi}

\begin{example}
If we take \textbf{\textit{Vect}}$_K$ with the usual tensor product, since we assume $\car(K)\neq 2$, we have that $\textbf{\textit{BICat}}(\textbf{\textit{Lie}}(\textbf{\textit{Vect}}_K))$ and $\textbf{\textit{BICat}}(\textbf{\textit{LieAlg}}_K)$ are isomorphic as categories.
\end{example}
	
With this, we can define what is a braiding on a categorical Lie object in $\mathcal{LM}_K$.

\begin{defi}
	Let $\mathcal{C}=$({\scriptsize\begin{tikzcd}C_1\arrow[d,"f_1"]\\ D_1\end{tikzcd}}, {\scriptsize\begin{tikzcd}C_0\arrow[d,"f_0"]\\ D_0\end{tikzcd}}$,s,t,e,k$) be a categorical Lie object in $\mathcal{LM}_K$.
	
	A \emph{braiding }on $\mathcal{C}$ is a triple $\tau=(\tau^{C_0,D_0},\tau^{D_0,C_0},\tau^{2})$ where:
	\begin{itemize}
		\item $\tau^2\colon D_0\times D_0\xrightarrow{}D_1$ is a $K$-bilinear map such that $(D_1,D_0,s,t,e,k,\tau_2)$ is a braided crossed module of Lie $K$-algebras,
		\item $\tau^{D_0,C_0}\colon D_0\times C_0\xrightarrow{}C_1$ and $\tau^{C_0,D_0}\colon C_0\times D_0\xrightarrow{}C_1$ are $K$-bilinear maps which, with $\tau^2$, verifies the following properties for $c\in C_0$, $d,d'\in D_0$, $x\in C_1$, $y\in D_1$:
		\begin{align*}
		f_1(\tau^{C_0,D_0}_{c,d})=\tau^2_{f_0(c),d}\ & \quad \text{and} \quad  \ f_1(\tau^{D_0,C_0}_{d,c})=\tau^2_{d,f_0(c)},\\
		\tau^{C_0,D_0}_{c,d}\colon c*^{C_{0}}_{D_0}d\xrightarrow{}-c*^{C_{0}}_{D_0}d \ &\quad \text{and}  \quad \ 	\tau^{D_0,C_0}_{d,c}\colon -c*^{C_{0}}_{D_0}d\xrightarrow{}c*^{C_{0}}_{D_0}d.
		\end{align*}
		The following diagrams are verified in the internal category:
		\begin{center}
		\begin{tikzcd}
		{s_1(x)*^{C_0}_{D_0}s_2(y)}\arrow[d,"{\tau^{C_0,D_0}_{s_1(x),s_2(y)}}"]\arrow[r,"{x*^{C_1}_{C_0}y}"]& {t_1(x)*^{C_0}_{D_0}t_2(y)}\arrow[d,"{\tau^{C_0,D_0}_{t_1(x),t_2(y)}}"]\\
		{-s_1(x)*^{C_0}_{D_0}s_2(y)}\arrow[r,"{-x*^{C_1}_{D_1}y}"]& {-t_1(x)*^{C_0}_{D_0}t_2(y)}.
		\end{tikzcd},
		\begin{tikzcd}
		{-s_1(x)*^{C_0}_{D_0}s_2(y)}\arrow[d,"{\tau^{D_0,C_0}_{s_2(y),s_1(x)}}"]\arrow[r,"{-x*^{C_1}_{C_0}y}"]& {-t_1(x)*^{C_0}_{D_0}t_2(y)}\arrow[d,"{\tau^{D_0,C_0}_{t_2(y),t_1(x)}}"]\\
		{s_1(x)*^{C_0}_{D_0}s_2(y)}\arrow[r,"{x*^{C_1}_{D_1}y}"]& {t_1(x)*^{C_0}_{D_0}t_2(y)}.
		\end{tikzcd},
		\end{center}
	And we have the following properties:
		\begin{align*}
		\tau_{c,[d,d']_{D_0}}^{C_0,D_0}&=\tau_{c*^{C_0}_{D_0} d,d'}^{C_0,D_0}-\tau_{c*^{C_0}_{D_0} d',d}^{C_0,D_0},\\
		\tau_{[d,d']_{D_{0}},c}^{D_0,C_0}&=-\tau_{d,c*^{C_0}_{D_0} d'}^{D_0,C_0}-\tau_{c*^{C_0}_{D_0} d,d'}^{C_0,D_0},\\
		\tau_{c,[d,d']_{D_0}}^{C_0,D_0}&=\tau_{c*^{C_0}_{D_0} d,d'}^{C_0,D_0}+\tau_{d,c*^{C_0}_{D_0} d'}^{D_0,C_0},\\
		\tau_{[d,d']_{D_0},c}^{D_0,C_0}&=-\tau_{d,c*^{C_0}_{D_0} d'}^{D_0,C_0}+\tau_{d',c*^{C_0}_{D_0} d}^{D_0,C_0}.
		\end{align*}
	\end{itemize}
	We will say that ({\scriptsize\begin{tikzcd}C_1\arrow[d,"f_1"]\\ D_1\end{tikzcd}}, {\scriptsize\begin{tikzcd}C_0\arrow[d,"f_0"]\\ D_0\end{tikzcd}}$,s,t,e,k,\tau$) is a \emph{braided categorical Lie object in $\mathcal{LM}_K$}.
\end{defi}

\begin{remark}
	A braiding, in fact, is a pair $\tau=(\tau_1,\tau_2)$ but, for simplicity, we denote $\tau_2(d,d')=\tau^2_{d,d'}$ and $\tau_1\colon (C_0\otimes D_0)\oplus (D_0\otimes C_0)\xrightarrow{}C_1$ by the expression $\tau_1((c\otimes d)+(d'\otimes c'))=\tau_{c,d}^{C_0,D_0}+\tau_{d',c'}^{D_0,C_0}$.
\end{remark}

\begin{defi}
	Let ({\scriptsize\begin{tikzcd}C_1\arrow[d,"f_1"]\\ D_1\end{tikzcd}},{\scriptsize\begin{tikzcd}C_0\arrow[d,"f_0"]\\D_0\end{tikzcd}}$,s,t,e,k,\tau$) and ({\scriptsize\begin{tikzcd}C_1'\arrow[d,"g_1"]\\ D_1'\end{tikzcd}},{\scriptsize\begin{tikzcd}C_0'\arrow[d,"g_0"]\\ D_0'\end{tikzcd}}$,s',t',e',k',\psi$) be braided categorical Lie objects in $\mathcal{LM}_K$. A \emph{braided internal functor between categorical Lie objects in $\mathcal{LM}_K$} is an internal functor $((F^1_1,F^0_1),(F^1_0,F^0_0))$ between the respective categorical Lie objects which verify;
	\begin{itemize}
		\item $(F^0_1,F^0_0)\colon (D_1,D_0,s_2,t_2,e_2,k_2,\tau_2)\xrightarrow{}(D'_1,D'_0,s'_2,t'_2,e'_2,k'_2,\psi_2)$ is a braided internal functor between categorical Lie $K$-algebras.	
		\item $F^1_1(\tau_{c,d}^{C_0,D_0})=\psi_{F^1_0(c),F^0_0(d)}^{C_0',D_0'}$ for $c\in C_0$, $d\in D_0$.
		\item $F^1_1(\tau_{d,c}^{D_0,C_0})=\psi_{F^0_0(d),F^1_0(c)}^{D'_0,C'_0}$ for $c\in C_0$, $d\in D_0$.
	\end{itemize}
\end{defi}

We want to introduce the braiding for the categorical Leibniz $K$-algebras with the previous idea.

For this, as in the case for braiding of crossed modules of Leibniz $K$-algebras we will use two $K$-bilinear maps $\tau,\psi\colon C_0\times C_0\xrightarrow{}C_1$. Taking now as braiding in the inclusion Lie object in $\mathcal{LM}_K$, $\bar{\tau}$, we define $\bar{\tau}_{a,\overline{b}}^{C_0,\Lie(C_0)}=\tau_{a,b}$, $\bar{\tau}_{\overline{a},b}^{\Lie(C_0),C_0}=-\psi_{b,a}$ and $\tau^2_{\overline{a},\overline{b}}=\overline{\tau_{a,b}}=\overline{-\psi_{b,a}}$ where, again, we introduce a quotient in $C_1$ whose elements we will denote as $\overline{x}$.

Doing that we obtain the following definition.
	
\begin{defi}
	A \emph{braiding} for the categorical Leibniz $K$-algebra $(C_1,C_0,s,t,e,k)$ is a pair $(\tau,\psi)$ of $K$-bilinear maps $\tau,\psi\colon C_0\times C_0\xrightarrow{}C_1$, $(a,b)\mapsto \tau_{a,b}$ and $(a,b)\mapsto \psi_{a,b}$, verifying:
	\begin{equation}\label{LeiT1}
	\tau_{a,b}\colon [a,b]\xrightarrow{}-[a,b] \quad \text{and} \quad \psi_{a,b}\colon [a,b]\xrightarrow{}-[a,b], \tag{LeibT1}
	\end{equation}
	\begin{equation}\label{LeiT2}
	\begin{tikzcd}
	{[s(x),s(y)]}\arrow[d,"{\tau_{s(x),s(y)}}"]\arrow[r,"{[x,y]}"]& {[t(x),t(y)]}\arrow[d,"{\tau_{t(x),t(y)}}"]\\
	{-[s(x),s(y)]}\arrow[r,"{-[x,y]}"]& {-[t(x),t(y)]},
	\end{tikzcd} \ \ \
	\begin{tikzcd}
	{[s(x),s(y)]}\arrow[d,"{\psi_{s(x),s(y)}}"]\arrow[r,"{[x,y]}"]& {[t(x),t(y)]}\arrow[d,"{\psi_{t(x),t(y)}}"]\\
	{-[s(x),s(y)]}\arrow[r,"{-[x,y]}"]& {-[t(x),t(y)]},
	\end{tikzcd}\tag{LeibT2}
	\end{equation}
	\begin{align}
	\tau_{a,[b,c]} &=\tau_{[a,b],c}-\tau_{[a,c],b} ,\tag{LeibT3}\label{LeiT3}\\
	\psi_{a,[b,c]} &=\tau_{[a,b],c}-\psi_{[a,c],b},\tag{LeibT4}\label{LeiT4}\\
	\tau_{a,[b,c]} &=\tau_{[a,b],c}-\psi_{[a,c],b},\tag{LeibT5}\label{LeiT5}\\
	\psi_{a,[b,c]} &=\psi_{[a,b],c}-\psi_{[a,c],b},\tag{LeibT6}\label{LeiT6}
	\end{align}
	for all $a,b,c\in C_0$, $x,y\in C_1$.
	
	If the definition is verified, then we say that $(C_1,C_0,s,t,e,k,(\tau,\psi))$ is a \emph{braided categorical Leibniz $K$-algebra.}
\end{defi}

\begin{defi}
	Let $(C_1,C_0,s,t,e,k,(\tau,\psi))$ and $(C'_1,C'_0,s',t',e',k',(\tau',\psi'))$ be two braided categorical Leibniz $K$-algebras.
	
	A \emph{braided internal functor} between two braided categorical Leibniz $K$-algebras is an internal functor $(C_1,C_0,s,t,e,k)\xrightarrow{(F_1,F_0)}(C'_1,C'_0,s',t',e',k')$ verifying:
	\begin{align}
	F_1(\tau_{a,b})&=\tau'_{F_0(a),F_0(b)},\tag{LeibHT1}\label{LHT1}\\
	F_1(\psi_{a,b})&=\psi'_{F_0(a),F_0(b)},\tag{LeibHT2}\label{LHT2}
	\end{align}
	for $a,b\in C_0$.
	
	We denote the category of braided categorical Leibniz $K$-algebras and braided internal functors between them as $\textbf{\textit{BICat}}(\textbf{\textit{LeibAlg}}_K)$.
\end{defi}

We want to know how see the braided categorical Lie $K$-algebras as a particular case of braided categorical Leibniz $K$-algebras. We will se that in the following property.

\begin{prop}
	Let $C_1$ and $C_0$ be Lie $K$-algebras. Then,
	$(C_1,C_0,s,t,e,k,\tau)$ is a braided categorical Lie $K$-algebra if and only if $(C_1,C_0,s,t,e,k,(\tau,\tau^{-}))$ is a braided categorical Leibniz $K$-algebra.
	
	$\tau^{-}\colon C_0\times C_0\xrightarrow{}C_1$ is defined as $\tau^{-}_{a,b}=-\tau_{b,a}$.
	\begin{proof}
		It is immediate using the anticommutativity that \eqref{LeiT1} and \eqref{LeiT2} are reduced respectively to \eqref{LieT1} and \eqref{LieT2}.
		
		For the four remaining properties we can se that they are equal to the last four of crossed modules, and the last two of braiding for categorical Lie $K$-algebras are equal to the last two of their respective crossed modules (inverted). For that we can argue in the same way we do with braided for crossed modules of Leibniz $K$-algebras as follows.
		
		It is clear that \eqref{LeiT3} and \eqref{LieT4} are identical, and is very easy to prove that \eqref{LeiT6} is equivalent to \eqref{LieT3}.

		To see the last equivalences we must prove an earlier property, which is verified for both braidings under our assumptions.
		
		If $a,b,c\in C_0$, then $\tau_{[a,b],c}=-\tau_{c,[a,b]}$.
		
		In the Lie case we have it using  Proposition~\ref{trenzaanticoncor}.
		
		In Leibniz case it is not true in general, because we need $\tau_{a,b}=-\psi_{b,a}$.
		
		Using \eqref{LeiT4} and \eqref{LeiT5} we can observe that
 \[
		\tau_{a,[b,c]}=\psi_{a,[b,c]}=\tau^{-}_{a,[b,c]}=-\tau_{[b,c],a},
\]
which is the property we need to prove.
		
		With this property we can prove easily the remaining equivalences. That is that \eqref{LeiT4} is equivalent to \eqref{LieT3}, and \eqref{LeiT5} is equivalent to \eqref{BLie4}.
	\end{proof}
\end{prop}

The next proposition gives the construction which yields us to obtain the definition.

\begin{prop}
	Let $(C_1,C_0,s,t,e,k,(\tau,\psi))$ be a braided categorical Leibniz $K$-algebra.
	
	{\rm ({\scriptsize\begin{tikzcd}C_1\arrow[d,"\pi_{C_1}"]\\\frac{C_1}{[\tau_{C_0,C_0}]}\end{tikzcd}},{\scriptsize\begin{tikzcd}C_0\arrow[d,"\pi_{C_0}"]\\\Lie(C_0)\end{tikzcd}}
$,(s,\overline{s}),(t,\overline{t}),(e,\overline{e}),(k,\widetilde{k}),\bar{\tau}$)} is a braided categorical Lie object in $\mathcal{LM}_K$,
 where we denote as $\frac{C_1}{[\tau_{C_0,C_0}]}$ the Lie $K$-algebra which is a Leibniz quotient of $C_1$ by the ideal generated by elements of the  form
  $[x,x]$ and $\tau_{a,b}+\psi_{b,a}$ with $x\in C_1$, $a,b\in C_0$, its elements as $\overline{x}$ with $x\in C_1$, and the maps are the following ones:
	\begin{itemize}
		\item $\overline{s}\colon \frac{C_1}{[\tau_{C_0,C_0}]}\xrightarrow{}\Lie(C_0)$ defined as $\overline{s}(\overline{x})=\overline{s(x)}$ for $\overline{x}\in \frac{C_1}{[\tau_{C_0,C_0}]}$;
		\item $\overline{t}\colon \frac{C_1}{[\tau_{C_0,C_0}]}\xrightarrow{}\Lie(C_0)$ defined as $\overline{t}(\overline{x})=\overline{t(x)}$ for $\overline{x}\in \frac{C_1}{[\tau_{C_0,C_0}]}$;
		\item $\overline{e}\colon \Lie(C_0)\xrightarrow{}\frac{C_1}{[\tau_{C_0,C_0}]}$ defined as $\overline{e}(\overline{a})=\overline{e(a)}$ for $\overline{a}\in \Lie(C_0)$;
		\item $\widetilde{k}\colon \frac{C_1}{[\tau_{C_0,C_0}]}\times_{\Lie(C_0)}\frac{C_1}{[\tau_{C_0,C_0}]} \xrightarrow{}\frac{C_1}{[\tau_{C_0,C_0}]}$ defined as $\overline{k}((\overline{x},\overline{y}))=\overline{\mathring{k}(x,y)}$ for $(\overline{x},\overline{y})\in \frac{C_1}{[\tau_{C_0,C_0}]}\times_{\Lie(C_0)}\frac{C_1}{[\tau_{C_0,C_0}]}$, where $\mathring{k}$ is again the extension to the product $\mathring{k}(x,y)=x+y-e(s(y))$ (we can take $\mathring{k}'(x,y)=x+y-e(t(x))$ too, because in the quotient it will not change anything);
		\item $\bar{\tau}^{C_0,\Lie(C_0)}\colon C_0\times \Lie(C_0)\xrightarrow{}C_1$ defined as $\bar{\tau}_{a,\overline{b}}^{C_0,\Lie(C_0)}=\tau_{a,b}$ for $a\in C_0$, $\overline{b}\in \Lie(C_0)$;
		\item $\bar{\tau}^{\Lie(C_0),C_0}\colon \Lie(C_0)\times C_0\xrightarrow{}C_1$ defined as $\bar{\tau}_{\overline{a},b}^{\Lie(C_0),C_0}=-\psi_{b,a}$ for $\overline{a}\in \Lie(C_0)$, $b\in C_0$;
		\item $\bar{\tau}^2\colon \Lie(C_0)\times \Lie(C_0)\xrightarrow{}\frac{C_1}{[\tau_{C_0,C_0}]}$ defined as $\tau^2_{\overline{a},\overline{b}}=\overline{\tau_{a,b}}=\overline{-\psi_{b,a}}$ for $\overline{a},\overline{b}\in \Lie(C_0)$.
	\end{itemize}
\end{prop}

\begin{remark}
	The bottom part $(\frac{C_1}{[\tau_{C_0,C_0}]},\Lie(C_0),\overline{s},\overline{t},\overline{e},\widetilde{k},\bar{\tau}_2)$ will be called Lieization, and it is again functorial.
	
	If we apply this Lieization on a braided categorical Lie $K$-algebra, thought as a crossed module of Leibniz $K$-algebras with the action with the braiding $(\tau,\tau^{-})$, the new generator is null
 \[
	\tau_{a,b}+\psi_{b,a}=\tau_{a,b}+\tau^{-}_{b,a}=\tau_{a,b}-\tau_{a,b}=0.
\]
	For that, again, in the Lie case, we obtain that doing the Lieization to and braided categorical Lie $K$-algebra gives us something trivially isomorphic to identity.
\end{remark}

\begin{prop}
	If {\rm ({\scriptsize\begin{tikzcd}C_1\arrow[d,"f_1"]\\D_1\end{tikzcd}}
	,{\scriptsize\begin{tikzcd}C_0\arrow[d,"f_0"]\\D_0\end{tikzcd}}$,s,t,e,k,\tau$)} is a braided categorical Lie object in $\mathcal{LM}_K$,
then $(C_1,C_0,s_1,t_1,e_1,k_1,(\bar{\tau}^\tau,\bar{\psi}^\tau))$ is a braided categorical Leibniz $K$-algebra, where $[x,y]_{C_1}=x*^{C_1}_{D_1}y$ and $[a,b]_{C_0}=a*^{C_0}_{D_0}b$ for $x,y\in C_1, a,b\in C_0$ and $\bar{\tau}^\tau_{a,b}=\tau^{C_0,D_0}_{a,f_0(b)}$, $\bar{\psi}^\tau_{a,b}=-\tau^{D_0,C_0}_{f(b),a}$ for $a,b\in C_0$.
\end{prop}

We have again a functor $BI\Phi\colon\textbf{\textit{BICat}}(\textbf{\textit{LeibAlg}}_K)\xrightarrow{}\textbf{\textit{BICat}}(\textbf{\textit{Lie}}(\mathcal{LM}_K))$ that is full, and another $BI\Psi\colon\textbf{\textit{BICat}}(\textbf{\textit{Lie}}(\mathcal{LM}_K))\xrightarrow{}\textbf{\textit{BICat}}(\textbf{\textit{LeibAlg}}_K)$, verifying $BI\Psi\circ BI\Phi=\Id_{\textbf{\textit{BICat}}(\textbf{\textit{LeibAlg}}_K)}$. $BI\Phi$ is a full inclusion functor.

\section{The equivalence between the categories of braided crossed modules and braided internal categories in the case of Leibniz algebras} \label{S:Leibequiv}

In the following sections we will prove other properties to verify that the two definitions given works properly.

The first one that we want to verify with our definition is that the categories $\textbf{\textit{BICat}}(\textbf{\textit{LeibAlg}}_K)$ and $\textbf{\textit{BX}}(\textbf{\textit{LeibAlg}}_K)$ must be equivalent, as in the case of groups and Lie $K$-algebras.

Also the equivalence must generalise the case of Lie $K$-algebras (this is, in the case in which we restrict the Leibniz $K$-algebras we want that the braidings verifies $\{n,n'\}=-\langle n',n\rangle$ and $\tau_{a,b}=-\psi_{b,a}$ and the functors for the Lie case would be recovered) and must be a extension of the one given to the non-braiding case.

With this in mind we will prove the equivalence. Remember that we are working with $\car(K)\neq 2$.

\begin{prop}
	Let $\mathcal{X}=(M,N,(\cdot_1,\cdot_2),\partial,(\{-,-\},\langle-,-\rangle))$ be a braided crossed module of Leibniz $K$-algebras.
	
	Then $\mathcal{C}_\mathcal{X}\coloneqq(M\rtimes N,N,\bar{s},\bar{t},\bar{e},\bar{k},(\bar{\tau},\bar{\psi}))$ is a braided categorical Leibniz $K$-algebra where:
	\begin{itemize}
		\item $\bar{s}\colon M\rtimes N\xrightarrow{} N$, \ $\bar{s}((m,n))=n$,
		\item $\bar{t}\colon M\rtimes N\xrightarrow{} N$, \ $\bar{t}((m,n))=\partial m+n$,
		\item $\bar{e}\colon N\xrightarrow{} M\rtimes N$,  \ $\bar{e}(n)=(0,n)$,
		\item $\bar{k}\colon (M\rtimes N)\times_N (M\rtimes N)\xrightarrow{}M\rtimes N$, where the source is the pullback of $\bar{t}$ with $\bar{s}$, defined as $k(((m,n),(m',\partial m+n)))=(m+m',n)$,
		\item $\bar{\tau}\colon N\times N\xrightarrow{}M\rtimes N$,  \  $\bar{\tau}_{n,n'}=(-2\{n,n'\},[n,n'])$,
		\item $\bar{\psi}\colon N\times N\xrightarrow{}M\rtimes N$, \  $\bar{\psi}_{n,n'}=(-2\langle n,n'\rangle,[n,n'])$.
	\end{itemize}
	\begin{proof}
		Other than the braiding, it has already been proven that $(M\rtimes N,N,\bar{s},\bar{t},\bar{e},\bar{k})$ is a categorical Leibniz $K$-algebra as can be seen  in \cite{ThRa}. We only need to check the braiding axioms for this internal category.
		
		We will start with \eqref{LeiT1}. Let  $n,n'\in N$.
		\begin{align*}
		\bar{s}(\bar{\tau}_{n,n'})&=\bar{s}((-2\{n,n'\},[n,n']))=[n,n'],\\
		\bar{t}(\bar{\tau}_{n,n'})&=\bar{t}((-2\{n,n'\},[n,n']))=-2\partial\{n,n'\}+[n,n']\\
		&=-2[n,n']+[n,n']=-[n,n'],
		\end{align*}
		where we use \eqref{LeiB1}. In the same way we can prove this property of $\bar{\psi}$ by the symmetry of the construction.
		
		We will prove now \eqref{LeiT2}. Again we will only check this for $\bar{\tau}$.
		 Let $x=(m,n),y=(m',n')\in M\rtimes N$.
		
		We need to show that $\tau_{t(x),t(y)}\circ [x,y]=-[x,y]\circ \tau_{s(x),s(y)}$. For that we will write the equalities in function of the data given by the braided crossed module.
\allowdisplaybreaks
		\begin{align*}
		&\tau_{t(x),t(y)}\circ [x,y]\\
		&=\bar{k}(([(m,n),(m',n')],(-2\{\bar{t}((m,n)),\bar{t}((m',n'))\},[\bar{t}((m,n)),\bar{t}((m',n'))])))\\
		&=\bar{k}(([(m,n),(m',n')],(-2\{\partial m+n,\partial m'+n'\},[\partial m+n,\partial m'+n'])))\\
		&=\bar{k}((([m,m']+n\cdot_1 m'+m\cdot_2 n',[n,n']),(-2\{\partial m+n,\partial m'+n'\},\\
		&\qquad\qquad\qquad\qquad\qquad\qquad\qquad\qquad\qquad\qquad[\partial m+n,\partial m'+n'])))\\
		&=([m,m']+n\cdot_1 m'+m\cdot_2n'-2\{\partial m+n,\partial m'+n'\},[n,n'])\\
		&=([m,m']+n\cdot_1 m'+m\cdot_2n'-2\{\partial m,\partial m'\}-2\{\partial m,n'\}-2\{n,\partial m'\}\\			&\qquad\qquad\qquad\qquad\qquad\ \quad\qquad\qquad\qquad\qquad\qquad-2\{n,n'\},[n,n'])\\
		&=([m,m']+n\cdot_1 m'+m\cdot_2n'-2[m,m']-2(m\cdot_2n')-2(n\cdot_1 m')\\
		&\qquad\qquad\qquad\qquad\qquad\qquad \ \quad\qquad\qquad\qquad\qquad-2\{n,n'\},[n,n'])\\
		&=(-[m,m']-n\cdot_1 m'-m\cdot_2n'-2\{n,n'\},[n,n']),
		\end{align*}
		where we use \eqref{LeiB2}, \eqref{LeiB3} and \eqref{LeiB4} in the sixth equality. In the other way,
		\begin{align*}
		& -[x,y]\circ \tau_{s(x),s(y)}\\
		&=\bar{k}(((-2\{\bar{s}((m,n)),\bar{s}((m,n'))\},[\bar{s}((m,n)),\bar{s}((m',n'))]),-[(m,n),(m',n')]))\\
		&=\bar{k}(((-2\{n,n'\},[n,n']),-[(m,n),(m',n')]))\\	
		&=\bar{k}(((-2\{n,n'\},[n,n']),(-[m,m']-n\cdot_1 m'-m\cdot_2 n',-[n,n'])))\\
		&=(-2\{n,n'\}-[m,m']-n\cdot_1 m'-m\cdot_2n',[n,n']).
		\end{align*}

		We will verify \eqref{LeiT3} below. Let $n,n',n''\in N$. Then
		\begin{align*}
		\bar{\tau}_{n,[n',n'']} & =(-2\{n,[n',n'']\},[n,[n',n'']])\\
		&=(-2(\lbrace[n,n'],n'' \rbrace-\lbrace [n,n''],n'\rbrace),[[n,n'],n'']-[[n,n''],n'])\\
		&=(-2\lbrace[n,n'],n'' \rbrace,[[n,n'],n''])-(-2\lbrace [n,n''],n'\rbrace,[[n,n''],n'])\\
		&=\bar{\tau}_{[n,n'],n''}-\bar{\tau}_{[n,n''],n'},
		\end{align*}
		where we use \eqref{LeiB5} and the Leibniz identity in the second equality.
		
		The same proof is true, for the symmetry of properties, using \eqref{LeiB8} for \eqref{LeiT6}.
		
		Finally, we will show that \eqref{LeiT4} and \eqref{LeiT5} are verified.
		\begin{align*}
		\bar{\psi}_{n,[n',n'']} & =(-2\langle n,[n',n'']\rangle,[n,[n',n'']])\\
		&=(-2(\{[n,n'],n''\}-\langle[n,n''],n'\rangle),[[n,n'],n'']-[[n,n''],n'])\\
		&=(-2\{[n,n'],n''\},[[n,n'],n''])-(-2\langle[n,n''],n'\rangle,[[n,n''],n'])\\
		&=\bar{\tau}_{[n,n'],n''}-\bar{\psi}_{[n,n''],n'}\\
		&=(-2(\{[n,n'],n''\}-\langle[n,n''],n'\rangle),[[n,n'],n'']-[[n,n''],n'])\\
		&=(-2\{ n,[n',n'']\},[n,[n',n'']])=\bar{\tau}_{n,[n',n'']},
		\end{align*}
		where we use \eqref{LeiB6} along with the Leibniz identity in the second equality; and \eqref{LeiB7} with the Leibniz identity in the penultimate equality.
		
		Thus the braiding axioms are verified for the categorical Leibniz $K$-algebra.
	\end{proof}
\end{prop}

\begin{remark}
	Note that if $\mathcal{X}$ is, in fact, a braided crossed module of Lie $K$-algebras, then
	\[
	\bar{\tau}_{n,n'}=(-2\{n,n'\},[n,n'])=-(-2\langle n',n\rangle,[n',n])=-\bar{\psi}_{n',n}
	\]
	and we recover the construction for the Lie case (see \cite{TFM}).
\end{remark}

\begin{prop}
	We have a functor $\mathcal{C}\colon\textbf{\textit{BX}}(\textbf{\textit{LeibAlg}}_K)\xrightarrow{}\textbf{\textit{BICat}}(\textbf{\textit{LeibAlg}}_K)$ defined as
	\[\mathcal{C}(\mathcal{X}\xrightarrow{(f_1,f_2)}\mathcal{X}') \coloneqq \mathcal{C}_\mathcal{X}\xrightarrow{(f_1\times f_2,f_2)}\mathcal{C}_{\mathcal{X}'},\]
	where $\mathcal{C}_\mathcal{X}$ is defined in the previous proposition.
	\begin{proof}
		We know that the pair $(f_1\times f_2,f_2)$ is an internal functor between the respective internal categories, since what we are trying to do is to extend an existing functor (see \cite{ThRa})
 to the braided case. In the same way, as it is the same functor, we already know that, if it is well defined, it verifies the properties of functor,
  since the composition and identity are the same as in the categories without braiding.
		
		Because of that, to conclude this proof, is enough to see that $(f_1\times f_2,f_2)$ is a braided internal functor of braided categorical Leibniz $K$-algebras.
		
		We will verify \eqref{LHT1}. Let $n,n' \in N$.
		\begin{align*}
		(f_1\times f_2)(\bar{\tau}_{n,n'}) & =(f_1\times f_2)((-2\{n,n'\},[n,n']))=(-2f_1(\{n,n'\}),f_2([n,n']))\\
		&=(-2\{f_2(n),f_2(n')\}',[f_2(n),f_2(n')])=\bar{\tau}'_{f_2(n),f_2(n')},
		\end{align*}
		where we use \eqref{LHB1} in the penultimate equality.
		
		Again, because of symmetry of the braiding's properties and the construction, the same proof is true, using \eqref{LHB2} to prove \eqref{LHT2}.
	\end{proof}
\end{prop}

\begin{prop}
	Let $\mathcal{C}=(C_1,C_0,s,t,e,k,(\tau,\psi))$ be a braided categorical Leibniz $K$-algebra.
	
	Then $\mathcal{X}_\mathcal{C}\coloneqq (\ker(s),C_0,({}^{e}\cdot,\cdot^e),\partial_t,(\{-,-\}_\tau,\langle-,-\rangle_\psi))$ is a braided crossed module of Leibniz $K$-algebras where
	\begin{itemize}
		\item ${}^{e}\cdot\colon C_0\times \ker(s)\xrightarrow{}\ker(s)$, \ $a\ {{}^{e}\cdot} \ x\coloneqq[e(a),x]$,
		\item $\cdot^e\colon \ker(s)\times C_0\xrightarrow{}\ker(s)$, \  $x\cdot^e a\coloneqq [x,e(a)]$,
		\item $\partial_t\coloneqq t|_{\ker(s)}$,
		\item $\{-,-\}_\tau\colon C_0\times C_0\xrightarrow{}\ker(s)$, \ $\{a,b\}_\tau\coloneqq \frac{e([a,b])-\tau_{a,b}}{2}$,
		\item $\langle-,-\rangle_\psi\colon C_0\times C_0\xrightarrow{}\ker(s)$, \ $\langle a,b\rangle_\psi\coloneqq \frac{e([a,b])-\psi_{a,b}}{2}$.
	\end{itemize}
	\begin{proof}
		It is proven in \cite{ThRa} that, under these hypotheses, $(\ker(s),C_0,(\cdot^e,{}^{e}\cdot),\partial_t)$ is a crossed module of Leibniz $K$-algebras.
		
		For that, it is enough to show that $(\{-,-\}_\tau,\langle-,-\rangle_\psi)$ is a braiding on that crossed module.

		First let us see that it is well defined because the image falls in $C_1$ which is not $\ker(s)$.
 We will check it only to $\{-,-\}_\tau$ since for $\langle-,-\rangle_\psi$ we will have a completely symmetric argument. Let $a,b\in C_0$, and using \eqref{LeiT1}, we have
		\begin{align*}
		s(\{a,b\}_\tau)=s\big(\frac{e([a,b])-\tau_{a,b}}{2}\big)=\frac{[a,b]-[a,b]}{2}=0.
		\end{align*}

		Since they are well defined, we can see if they verify the properties. To check \eqref{LeiB1}--\eqref{LeiB4} we will only prove it for $\{-,-\}_\tau$,
 again because of the two braidings need to verify the same properties \eqref{LeiT1} and \eqref{LeiT2}, which are the ones we will use.
		
		To start we will check \eqref{LeiB1}. Let $a,b\in C_0$, and using \eqref{LeiT1}, we get
		\begin{align*}
		\partial_{t}\{a,b\}_\tau=t\big(\frac{e([a,b])-\tau_{a,b}}{2}\big)=\frac{[a,b]-(-[a,b])}{2}=\frac{2[a,b]}{2}=[a,b].
		\end{align*}

		We will see if it is verified \eqref{LeiB2}. Let $x,y\in \ker(s)$. Then
		\begin{align*}
		\{\partial_tx,\partial_ty\}_\tau=\frac{e([\partial_tx,\partial_ty])-\tau_{\partial_t x,\partial_t y}}{2}
		=\frac{e([t(x),t(y)])-\tau_{t(x),t(y)}}{2}.
		\end{align*}
		Let us see that $\frac{e([t(x),t(y)])-\tau_{t (x),t (y)}}{2}=[x,y]$.
		
		By the axiom \eqref{LeiT2} we know that we have the following equality:
		\[k(([x,y],\tau_{t(x),t(y)}))=k((\tau_{s(x),s(y)},-[x,y])).\]
		As $x\in \ker(s)$, we have that $s(x)=0$ (in the same way $y$), and $\tau_{s(x),s(y)}=0$ by $K$-bilinearity. We have then that
		\[k((\tau_{s(x),s(y)},-[x,y]))=k((0,-[x,y])),\]
		and therefore the following equality:
		\[k(([x,y],\tau_{t(x),t(y)}))=k((0,-[x,y])).\]
		Using now the $K$-linearity of $k$ in the previous expression, we obtain:
		\[0=k(([x,y],\tau_{t(x),t(y)}+[x,y])).\]
		
		Since $t(\tau_{t(x),t(y)}+[x,y])=-[t(x),t(y)]+[t(x),t(y)]=0=s(e(0))$ we can talk about $k((\tau_{t(x),t(y)}+[x,y],e(0)))$.
		Further $k((\tau_{t(x),t(y)}+[x,y],e(0)))=\tau_{t(x),t(y)}+[x,y]$ by the internal category axioms.
		
		Adding both equalities and by using the $K$-linearity of $k$ we get the following
		\[k(([x,y]+\tau_{t(x),t(y)}+[x,y],\tau_{t(x),t(y)}+[x,y]))=\tau_{t(x),t(y)}+[x,y].\]
		Therefore, by grouping, we have
		\[k((2[x,y]+\tau_{t(x),t(y)},\tau_{t(x),t(y)}+[x,y]))=\tau_{t(x),t(y)}+[x,y].\]
		By using that $\ker(s)$ is an ideal and the fact that $x$ or $y$ are in $\ker(s)$, we have $s(\tau_{t(x),t(y)}+[x,y])=[t(x),t(y)]-0=[t(x),t(y)]$, and so
 it makes sense to speak about the composition $k((e([t(x),t(y)]),\tau_{t(x),t(y)}+[x,y]))$, which is equal to $\tau_{t(x),t(y)}+[x,y]$.
		
		Subtracting both equalities and using the $K$-linearity of $k$, we obtain
		\[k((2[x,y]+\tau_{t(x),t(y)}-e([t(x),t(y)]),0))=0.\]
		
		Again, using the properties for internal categories, we have
		\begin{align*}
		0&=k((2[x,y]+\tau_{t(x),t(y)}-e([t(x),t(y)]),0))\\&=k((2[x,y]+\tau_{t(x),t(y)}-e([t(x),t(y)]),e(0)))
		\\&=2[x,y]+\tau_{t(x),t(y)}-e([t(x),t(y)]),
		\end{align*}
		which, since $\car(K)\neq 2$, gives us the required equality.

		As an observation to the above, in the part of the proof where we use that $x,y\in\ker(s) $, it is sufficient that one of the two is in that kernel.
		Therefore, by repeating the proof using this, we have the following equalities for $x \in \ker(s) $ and $y \in C_1 $:
		\begin{align*}
		\frac{e([t(x),t(y)])-\tau_{t (x),t (y)}}{2}=[x,y], \quad\frac{e([t(y),t(x)])-\tau_{t (y),t (x)}}{2}=[y,x].
		\end{align*}
		
		With these equalities, we will prove \eqref{LeiB3} and \eqref{LeiB4}.
		
		Let $a\in C_0$ and $x\in \ker(s)$. Then
		\begin{align*}
		&\{\partial_t x, a\}_\tau=\frac{e([t(x),t(e(a))])-\tau_{t (x),t(e(a))}}{2}=[x,e(a)]=x\cdot^e a,\\
		&\{ a,\partial_t x\}_\tau=\frac{e([t(e(a)),t(x)])-\tau_{t (e(a)),t(x)}}{2}=[e(a),x]=a\ {{}^{e}\cdot} \ x.
		\end{align*}
		
		We will see now the last conditions, starting with \eqref{LeiB5}. Let $a,b,c\in C_0$.
		\begin{align*}
		\{a,[b,c]\}_\tau & =\frac{e([a,[b,c]])-\tau_{a,[b,c]}}{2}=\frac{e([[a,b],c])-e([[a,c],b])-\tau_{[a,b],c}+\tau_{[a,c],b}}{2}\\
		&=\frac{e([[a,b],c])-\tau_{[a,b],c}}{2}-\frac{e([[a,c],b])-\tau_{[a,c],b}}{2}=\{[a,b],c\}_\tau-\{[a,c],b\}_\tau,
		\end{align*}
		where we use \eqref{LeiT3} and the Leibniz identity in the second equality.
		By symmetry we can prove \eqref{LeiB8}, using \eqref{LeiT6}.
		
		To conclude we will check \eqref{LeiB6} and \eqref{LeiB7}.
		\begin{align*}
		\langle a,[b,c]\rangle_\psi &=\frac{e([a,[b,c]])-\psi_{a,[b,c]}}{2}=\frac{e([[a,b],c])-e([[a,c],b])-\tau_{[a,b],c}+\psi_{[a,c],b}}{2}\\
		&=\frac{e([[a,b],c])-\tau_{[a,b],c}}{2}-\frac{e([[a,c],b])-\psi_{[a,c],b}}{2}=\{[a,b],c\}_\tau-\langle[a,c],b\rangle_\psi\\
		&=\frac{e([[a,b],c])-e([[a,c],b])-\tau_{[a,b],c}+\psi_{[a,c],b}}{2}=\frac{e([a,[b,c]])-\tau_{a,[b,c]}}{2}\\
		&=\{ a,[b,c]\}_\tau,
		\end{align*}
		where we use \eqref{LeiT4} in the second equality together with Leibniz identity and \eqref{LeiT5} in the penultimate equality with the Leibniz identity.
	\end{proof}
\end{prop}

\begin{remark}
	Note that if $\mathcal{C}$ is actually a braided categorical Lie $K$-algebra, then
	\[
	\{a,b\}_\tau=\frac{e([a,b])-\tau_{a,b}}{2}=-\frac{e([b,a])-\psi_{b,a}}{2}=-\langle b,a\rangle_\psi
	\]
	and we recover the construction for the Lie case.
\end{remark}

\begin{prop}
	We have a functor $\mathcal{X}\colon \textbf{\textit{ICat}}(\textbf{\textit{LeibAlg}}_K)\xrightarrow{}\textbf{\textit{BX}}(\textbf{\textit{LeibAlg}}_K)$ defined as
	\[\mathcal{X}(\mathcal{C}\xrightarrow{(F_1,F_0)}\mathcal{C}')=\mathcal{X}_\mathcal{C}\xrightarrow{(F_1^s,F_0)}\mathcal{X}_{\mathcal{C}'},\]
	where $\mathcal{X}_\mathcal{C}$ is defined in the previous proposition and $F_1^s\colon \ker(s)\xrightarrow{}\ker(s')$ is defined as $F_1^s(x)=F_1(x)$ for $x\in  \ker(s)$.
	\begin{proof}
		The fact that it is a functor between the categories without braiding is already shown \cite{ThRa},
 so we have to see that it can be extended to the braided case. For this, we have to verify the axioms of the homomorphisms of braided crossed of Leibniz $K$-algebras.
		
		We will start with \eqref{LHB1}. Let $a,b \in C_0$.
		\begin{align*}
		 F^s_1(\{a,b\}_\tau)&=F_1\big(\frac{e([a,b])-\tau_{a,b}}{2}\big)=\frac{F_1(e([a,b]))-F_1(\tau_{a,b})}{2}\\
		&=\frac{e'(F_0([a,b]))-\tau'_{F_0(a),F_0(b)}}{2}=\frac{e'([F_0(a),F_0(b)])-\tau'_{F_0(a),F_0(b)}}{2}\\
		&=\{F_0(a),F_0(b)\}_{\tau'},
		\end{align*}
		where we use \eqref{LHT1} in the third equality.
		
		Again by symmetry using \eqref{LHT2} we can prove using the same argument \eqref{LHB2}.
	\end{proof}
\end{prop}

\begin{remark}
	Note that, if $(M,N,(\cdot_1,\cdot_2),\partial,(\{-,-\},\langle-,-\rangle))$ then $\ker(\bar{s})=\{(m,0)\in M\rtimes N\mid m\in M\}$, where $\bar{s}$ is defined for the functor $\mathcal{C}$.
\end{remark}

\begin{prop}
	The categories $\textbf{\textit{BX}}(\textbf{\textit{LeibAlg}}_K)$ and $\textbf{\textit{ICat}}(\textbf{\textit{LeibAlg}}_K)$ are equivalent categories.
	
	Further, the functors $\mathcal{C}$ and $\mathcal{X}$ are inverse equivalences, where the natural isomorphisms
$\Id_{\textbf{\textit{BX}}(\textbf{\textit{LeibAlg}}_K)}\stackrel{\alpha}{\cong}\mathcal{X}\circ \mathcal{C}$ and $\Id_{\textbf{\textit{ICat}}(\textbf{\textit{LeibAlg}}_K)}\stackrel{\beta}{\cong}\mathcal{C}\circ \mathcal{X}$ are given by:
	\begin{itemize}
   \item If $\mathcal{Z}=(M,N,(\cdot_1,\cdot_2),\partial,(\{-,-\},\langle-,-\rangle))$ is a braided crossed module of Leibniz $K$-algebras, then $\alpha_{\mathcal{Z}}=(\alpha_M,\Id_N)$, where $\alpha_M\colon M\xrightarrow{} (M,0)$ is defined by $\alpha_M(m)=(m,0)$;
   \item If $\mathcal{D}=(C_1,C_0,s,t,e,k,(\tau,\psi))$ is a braided categorical Leibniz $K$-algebra, then $\beta_{\mathcal{D}}=(\beta_{s},\Id_{C_0})$, where $\beta_{C_1}\colon C_1\xrightarrow{}\ker(s)\rtimes C_0$ is defined by $\beta_{C_1}(x)=(x-e(s(x)),s(x))$.
 \end{itemize}	
	\begin{proof}
		It can be seen in \cite{ThRa} that they are well defined maps and that they are isomorphisms in the categories without braiding, as well as they are natural isomorphisms.
		
		For that, is enough to show that they are isomorphisms between braided objects.
		
		Immediately from definition, as in the crossed module and categorical case, they are isomorphisms if they are bijective morphisms,
since the inverse map verify the braided axioms for morphisms in their respectively category.
		
		We know that they are bijective maps, since they are isomorphisms between the categories without braiding. For that, we only have to verify that they are, in fact, morphisms.
		
		Let $\mathcal{Z}=(M,N,(\cdot_1,\cdot_2),\partial,(\{-,-\},\langle-,-\rangle))$ a braided crossed module of Leibniz $K$-algebras.
 Let's see that $\alpha_{\mathcal{Z}}=(\alpha_M,\Id_N)$ verifies \eqref{LHB1}. For that, we will take $n,n'\in N$:
		\begin{align*}
		\Id_N(\{n,n'\}_{\bar{\tau}})&=\{n,n'\}_{\bar{\tau}}=\frac{\bar{e}([n,n'])-\bar{\tau}_{n,n'}}{2}=\frac{(0,[n,n'])-(-2\{n,n'\},[n,n'])}{2}\\
		&=\frac{(2\{n,n'\},0)}{2}=(\{n,n'\},0)=\alpha_{M}(\{n,n'\}).
		\end{align*}
		
		Analogously \eqref{LHB2} is proven, by the similarity of definitions.
		
		Let $\mathcal{D}=(C_1,C_0,s,t,e,k,(\tau,\psi))$ be a braided categorical Leibniz $K$-algebra.
 We will check that $\beta_{\mathcal{D}}=(\beta_{s},\Id_{C_0})$ verifies \eqref{LHT1} and \eqref{LHT2}. For that, we only show the proof for \eqref{LHT1}, since the one for \eqref{LHT2} is identical.

		Let us consider $a,b\in C_0$. We have:
		\begin{align*}
		\Id_{C_0}(\bar{\tau}_{a,b})&=\bar{\tau}_{a,b}=(-2\{a,b\}_\tau,[a,b])=(-2\frac{e([a,b])-\tau_{a,b}}{2},[a,b])\\
		&=(\tau_{a,b}-e([a,b]),[a,b])=(\tau_{a,b}-e(s(\tau_{a,b})),s(\tau_{a,b}))=\beta_{C_1}(\tau_{a,b}).
		\end{align*}

		Therefore, since they are morphisms, we know that these are natural isomorphisms as we explained above and the equivalence of categories is obtained.
	\end{proof}
\end{prop}

\section{The non-abelian tensor product as example of braiding}\label{S:nonabeliantensorLeib}

As a last point that the definition verifies reasonable conditions, it is immediate to check that if $(M,[-,-])$ is a Leibniz $K$-algebra,
 then $([-,-],[-,-])$ is a braiding on $(M,M,([-,-],[-,-]),\Id_M)$.
This example is the analogous for the case of Leibniz $K$-algebras of the examples $(G,G,\Conj,\Id_G,[-,-])$ for groups and $(M,M,[-,-],\Id_M,[-,-])$ for Lie $K$-algebras.
 Further, this example generalizes the Lie example, since $[y,x]=-[x,y]$ in the Lie case.

We will give another symmetric example in the three constructions: the non-abelian tensor product. The non-abelian tensor product of Leibniz $K$-algebras was introduced by Gnedbaye in \cite{Gned}.
 In this paper the tensor product is denoted as $M \star N$, notation that we will keep, and its generators as $m*n $ and  $n*m$.
  In the general case it does not give rise to confusion, but in the case that we are interested is $M=N$ and there  generators would be denoted in the same way, giving rise to confusion. To avoid this we change the nomenclature, denoting $m*n$ as $m\otimes n$ and $n*m$ as $n \circledast m$, which differentiate the generators in any case.

\begin{defi}
	Let $M$ and $N$ two Leibniz $K$-algebras together with two Leibniz actions $\cdot=(\cdot_1,\cdot_2)$ of $M$ on $N$ and $*=(*_1,*_2)$ of $N$ on $M$.
	
	The \emph{non-abelian tensor product of $M$ and $N$}, denoted by $M\star N$, is the Leibniz $K$-algebra generated by the symbols $m\otimes n$ and $n\circledast m$ with $m\in M$, $n\in N$, together with the relations:
	\begin{align*}\label{RTLei1}
	\lambda(m\otimes n)&=\lambda m\otimes n=m\otimes \lambda n, \tag{RTLeib1}\\
	\lambda(n\circledast m)&=\lambda n\circledast m=n\circledast \lambda m,
	\end{align*}
	\begin{align*}\label{RTLei2}
	(m+m')\otimes n&=m\otimes n+m'\otimes n, \tag{RTLeib2}\\
	m\otimes (n+n')&=m\otimes n+m\otimes n',\\
	(n+n')\circledast m&=n\circledast m+n'\circledast m,\\
	n\circledast (m+m')&=n\circledast m+n\circledast m',
	\end{align*}	
	\begin{align*}\label{RTLei3}
	m\otimes [n,n']&=(m *_2 n)\otimes n'-(m*_2n')\otimes n, \tag{RTLeib3}\\
	n\circledast [m,m']&=(n \cdot_2 m)\circledast m'-(n\cdot_2m')\circledast m,\\
	[m,m']\otimes n&=(m\cdot_1 n)\circledast m'-m\otimes (n\cdot_2 m'),\\
	[n,n']\circledast m&=(n*_1 m)\otimes n'-n\circledast (m*_2 n'),
	\end{align*}
	\begin{align*}\label{RTLei4}
	m\otimes (m'\cdot_1 n)&=-m\otimes (n\cdot_2 m'),\tag{RTLeib4}\\
	n\circledast (n'*_1 m)&=-n\circledast (m*_2 n'),
	\end{align*}
	\begin{align*}\label{RTLei5}
	(m*_2 n)\otimes (m'\cdot_1 n')&=[m\otimes n,m'\otimes n']=(m\cdot_1n)\circledast(m'*_2 n'), \tag{RTLeib5}\\
	(m*_2 n)\otimes (n'\cdot_2 m')&=[m\otimes n,n'\circledast m']=(m\cdot_1n)\circledast(n'*_1 m'),\\
	(n*_1 m)\otimes (n'\cdot_2 m')&=[n\circledast m,n'\circledast m']=(n\cdot_2m)\circledast(n'*_1 m'),\\
	(n*_1 m)\otimes (m'\cdot_1 n')&=[n\circledast m,m'\otimes n']=(n\cdot_2m)\circledast(m'*_2 n'),
	\end{align*}
	where $m,m'\in M$, $n,n'\in N$.
\end{defi}

Again we have in natural way a crossed module for this case, as can be seen in \cite{Gned}. We will only show a particular case, which is the interesting one for our example.

\begin{prop}[\cite{Gned}]
	Let $M$ be a Leibniz $K$-algebra.
	
	Then $(M\star M,M,(\cdot_1,\cdot_2),\partial)$ is a crossed module of Leibniz $K$-algebras, where $M\star M$ is the non-abelian tensor product of $M$ with itself using the actions given by the Leibniz bracket, with

 \begin{itemize}
   \item  the left action over generators is given by $m\cdot_1 (m_1\otimes m_2)=[m,m_1]\otimes m_2-[m,m_2]\circledast m_1$, $m\cdot_1 (m_1\circledast m_2)=[m,m_1]\circledast m_2-[m,m_2]\otimes m_1$;
   \item  the right action is given by the expressions $(m_1\otimes m_2)\cdot_2 m=[m_1,m]\otimes m_2+m_1\otimes [m_2,m]$, $(m_1\circledast m_2)\cdot_2 m=[m_1,m]\circledast m_2+m_1\circledast [m_2,m]$;
   \item the map $\partial$ defined over generators as $\partial(m_1\otimes m_2)=[m_1,m_2]=\partial(m_1\circledast m_2)$.
 \end{itemize}

\end{prop}

\begin{remark}
	We will show how are the relations \eqref{RTLei3}--\eqref{RTLei5} for the non-abelian tensor product $M\star M$ with the action $([-,-],[-,-])$ on itself:
	\begin{align*}
     m_1\otimes [m_2,m_3] &=[m_1,m_2]\otimes m_3-[m_1,m_3]\otimes m_2, && \text{(RTLeib3)} \\
	m_1\circledast [m_2,m_3] &=[m_1,m_2]\circledast m_3-[m_1,m_3]\circledast m_2,\\
	[m_1,m_2]\otimes m_3 &=[m_1,m_3]\circledast m_2-m_1\otimes [m_3,m_2],\\
	[m_1,m_2]\circledast m_3 &=[m_1,m_3]\otimes m_2-m_1\circledast [m_3,m_2],\\
	m_1\otimes [m_2,m_3] & =-m_1\otimes [m_3,m_2], && \text{(RTLeib4)}\\
	m_1\circledast [m_2,m_3] & =-m_1\circledast [m_3, m_2], \\
	[m_1,m_2]\otimes [m_3,m_4] &=[m_1\otimes m_2,m_3\otimes m_4]=[m_1,m_2]\circledast[m_3,m_4], && \text{(RTLeib5)}\\
	[m_1,m_2]\otimes [m_3,m_4] &=[m_1\otimes m_2,m_3\circledast m_4]=[m_1,m_2]\circledast[m_3,m_4],\\
	[m_1,m_2]\otimes [m_3,m_4] &=[m_1\circledast m_2,m_3\circledast m_4]=[m_1,m_2]\circledast[m_3,m_4],\\
	[m_1,m_2]\otimes [m_3,m_4] &=[m_1\circledast m_2,m_3\otimes m_4]=[m_1,m_2]\circledast[m_3,m_4],
	\end{align*}
	where $m_1,m_2,m_3,m_4\in M$.
\end{remark}

The following example shows the necessity of a pair of braidings for the Leibniz $K$-algebras case since they will be different; and it is, how it was expected, analogous to the previous cases.

\begin{example}
	Let $M$ be a Leibniz $K$-algebra.
	
	The pair of $K$-bilinear maps $\{-,-\},\langle-,-\rangle\colon M\times M\xrightarrow{} M\otimes M$ defined as $\{m_1,m_2\}=m_1\otimes m_2$ and $\langle m_1,m_2\rangle=m_1\circledast m_2$ is a braiding on the crossed module of Leibniz $K$-algebras $(M\star M,M,(\cdot_1,\cdot_2),\partial)$, being this the one defined in the previous proposition.
	
	Let us check this. We will start with \eqref{LeiB1} taking $m_1,m_2\in M$,
	\[
	\partial \{m_1,m_2\}=\partial(m_1\otimes m_2)=[m_1,m_2]=\partial(m_1\circledast m_2)=\partial\langle m_1,m_2\rangle.
	\]
	
	The following properties will be proven on generators of $M\star M$, since the $K$-linearity and $K$-bilinearity of the maps let us say that if the properties are verified by generators, then are verified for the general elements (again, as in the Lie case, the bracket can be expressed as a sum of generators).
	
	With this we will prove \eqref{LeiB2}.
	\begin{align*}
	\{\partial(m_1\otimes m_2),\partial(m_3\otimes m_4)\}&=[m_1,m_2]\otimes [m_3,m_4]=[m_1\otimes m_2,m_3\otimes m_4],\\
	\{\partial(m_1\otimes m_2),\partial(m_3\circledast m_4)\}&=[m_1,m_2]\otimes [m_3,m_4]=[m_1\otimes m_2,m_3\circledast m_4],\\
	\{\partial(m_1\circledast m_2),\partial(m_3\otimes m_4)\}&=[m_1,m_2]\otimes [m_3,m_4]=[m_1\circledast m_2,m_3\otimes m_4],\\
	\{\partial(m_1\circledast m_2),\partial(m_3\circledast m_4)\}&=[m_1,m_2]\otimes [m_3,m_4]=[m_1\circledast m_2,m_3\circledast m_4],\\
	\langle\partial(m_1\otimes m_2),\partial(m_3\otimes m_4)\rangle&=[m_1,m_2]\circledast [m_3,m_4]=[m_1\otimes m_2,m_3\otimes m_4],\\
	\langle\partial(m_1\otimes m_2),\partial(m_3\circledast m_4)\rangle&=[m_1,m_2]\circledast [m_3,m_4]=[m_1\otimes m_2,m_3\circledast m_4],\\
	\langle\partial(m_1\circledast m_2),\partial(m_3\otimes m_4)\rangle&=[m_1,m_2]\circledast [m_3,m_4]=[m_1\circledast m_2,m_3\otimes m_4],\\
	\langle\partial(m_1\circledast m_2),\partial(m_3\circledast m_4)\rangle&=[m_1,m_2]\circledast [m_3,m_4]=[m_1\circledast m_2,m_3\circledast m_4],
	\end{align*}
	where in all cases we used \eqref{RTLei5}.
	
	Before to check the following axioms, we need to check a property that can be proven using \eqref{RTLei3} and \eqref{RTLei4}.
	
	Using \eqref{RTLei4} in the last equality of relation \eqref{RTLei3} and rewriting that equality and the second one, we get
	\begin{align*}
	m_1\circledast [m_2,m_3]=[m_1,m_2]\circledast m_3-[m_1,m_3]\circledast m_2,\\
	m_1\circledast [m_2,m_3]=[m_1,m_2]\circledast m_3-[m_1,m_3]\otimes m_2.
	\end{align*}
	Subtracting, we obtain the equality $[m_1,m_3]\otimes m_2=[m_1,m_3]\circledast m_2$.
	Using this last equality and the first and second equality of \eqref{RTLei3},  we obtain
	\begin{align*}
	m_1\otimes [m_2,m_3]&=[m_1,m_2]\otimes m_3-[m_1,m_3]\otimes m_2\\
	&=[m_1,m_2]\circledast m_3-[m_1,m_3]\circledast m_2=m_1\circledast [m_2,m_3].
	\end{align*}
	
	Let us verify now the first equality of \eqref{LeiB3} with $m,m_1,m_2\in M$,
	\begin{align*}
	\{\partial (m_1\otimes m_2),m \}&=[m_1,m_2]\otimes m=[m_1,m]\circledast m_2-m_1\otimes [m,m_2]\\
	&=[m_1,m]\circledast m_2+m_1\otimes [m_2,m]=[m_1,m]\otimes m_2+m_1\otimes [m_2,m]\\
	&=(m_1\otimes m_2)\cdot_2 m,
	\end{align*}
	where we use \eqref{RTLei3} and \eqref{RTLei4}.
	
	The second equality is analogous:
	\begin{align*}
	\{\partial (m_1\circledast m_2),m \}&=[m_1,m_2]\otimes m=[m_1,m]\circledast m_2+m_1\otimes [m_2,m]\\
	&=[m_1,m]\circledast m_2+m_1\circledast [m_2,m]=(m_1\circledast m_2)\cdot_2 m.
	\end{align*}
	Using again the exchange properties between $\otimes$ and $\circledast$, we will see the remaining equalities:
	\begin{align*}
	&\langle\partial (m_1\otimes m_2),m \rangle=[m_1,m_2]\circledast m=[m_1,m_2]\otimes m=(m_1\otimes m_2)\cdot_2 m,\\
	&\langle\partial (m_1\circledast m_2),m \rangle=[m_1,m_2]\circledast m=[m_1,m_2]\otimes m=(m_1\circledast m_2)\cdot_2 m.
	\end{align*}
	
	Now we will check the next axiom, \eqref{LeiB4}, where we will use again that we can exchange the symbols if in one side is the bracket. Starting with the first equality, we have
	\begin{align*}
	\{m,\partial (m_1\otimes m_2)\}&=m\otimes [m_1,m_2]=[m,m_1]\otimes m_2-[m,m_2]\otimes m_1\\
	&=[m,m_1]\otimes m_2-[m,m_2]\circledast m_1=m\cdot_1 (m_1\otimes m_2),
	\end{align*}
	where we use \eqref{RTLei3}. In an analogous way we obtain the second equality:
	\begin{align*}
	\{m,\partial (m_1\circledast m_2)\}&=m\otimes [m_1,m_2]=[m,m_1]\otimes m_2-[m,m_2]\otimes m_1\\
	&=[m,m_1]\circledast m_2-[m,m_2]\otimes m_1=m\cdot_1 (m_1\circledast m_2).
	\end{align*}
	With this the following properties are immediate:
	\begin{align*}
	&\langle m,\partial (m_1\otimes m_2)\rangle=m\circledast [m_1,m_2]=m\otimes [m_1,m_2]=m\cdot_1 (m_1\otimes m_2),\\
	&\langle m,\partial (m_1\circledast m_2)\rangle=m\circledast [m_1,m_2]=m\otimes [m_1,m_2]=m\cdot_1 (m_1\circledast m_2).
	\end{align*}
	
	To finalize we will prove \eqref{LeiB5}, because, if it is verified; \eqref{LeiB6}--\eqref{LeiB8} will be verified using the following properties:
	\begin{align*}
	&\{m,[m',m'']\}=m\otimes [m',m'']=m\circledast [m',m'']=\langle m,[m',m'']\rangle,\\
	&\{[m,m'],m''\}=[m,m']\otimes m''=[m,m']\circledast m''=\langle [m,m'],m'' \rangle.
	\end{align*}
	We have \eqref{LeiB5}:
	\begin{align*}
	\{m,[m',m'']\}=m\otimes [m',m'']&=[m,m']\otimes m''-[m,m'']\otimes m'\\
	&=\{[m,m'],m'' \}-\{[m,m''],m'\},
	\end{align*}
	where we used \eqref{RTLei3}.
\end{example}

\begin{remark}
	Note that the actions are reduced to a simplest notation, given by
	\begin{align*}
	& m\cdot_1(m_1\otimes m_2)=m\cdot_1(m_1\circledast m_2)=m\otimes [m_1,m_2]=m\circledast [m_1,m_2],\\
	& (m_1\otimes m_2)\cdot_2 m=(m_1\circledast m_2)\cdot_2 m=[m_1,m_2]\otimes m=[m_1,m_2]\circledast m.
	\end{align*}
\end{remark}

\begin{remark}
	This example generalizes the Lie example, since if we have that $m_1\circledast m_2=-m_2\otimes m_1$  as a new relation, we obtain the Lie tensor product of $M$ with itself using the adjoint action.
\end{remark}


\end{document}